\documentclass[a4paper,11pt]{article}

\usepackage[latin1]{inputenc}

\usepackage{latexsym}

\usepackage{amssymb}
\usepackage{amsfonts}
\usepackage{amsmath}
\usepackage{amsthm}
\usepackage{textcomp}

\title{Amoebas, tropical varieties and \\ compactification of Teichm\"uller spaces}
\author{Daniele Alessandrini \\ \small \ \\ \small \textit{Scuola Normale Superiore, Pisa, Italy} \\ \small \textit{E-mail address:} d.alessandrini@sns.it  }
\date{}

\newcommand{\nuovo}[1]{{{\bfseries \upshape #1}}}

\newcommand{\enne}{{\mathbb{N}}}
\newcommand{\ze}{{\mathbb{Z}}}
\newcommand{\qu}{{\mathbb{Q}}}
\newcommand{\erre}{{\mathbb{R}}}
\newcommand{\ci}{{\mathbb{C}}}

\newcommand{\cp}{{\mathbb{CP}}}

\newcommand{\cappa}{{\mathbb{K}}}
\newcommand{\effe}{{\mathbb{F}}}

\newcommand{\ident}{{\mbox{Id}}}%{{\mathbb{I}}}

\newcommand{\sld}{{SL_2(\ci)}}

\newcommand{\ideal}{{\mathfrak{I}}}

\newcommand{\ocors}{{\mathcal{O}}}

\newcommand{\plac}{{\mathcal{P}}}

\newcommand{\famil}{{\mathcal{F}}}
\newcommand{\gfamil}{{\mathcal{G}}}
\newcommand{\class}{{\mathcal{C}}}

\newcommand{\ameba}{{\mathcal{A}}}
\newcommand{\sferic}{{\mathcal{S}}}

\newcommand{\interno}[1]{{\stackrel{\circ}{#1}}}

\DeclareMathOperator{\tr}{tr}

\DeclareMathOperator{\homom}{Hom}

\DeclareMathOperator{\card}{Card}
\DeclareMathOperator{\Log}{Log}

\DeclareMathOperator{\segno}{sign}
\DeclareMathOperator{\rank}{rank}

\newenvironment{gmatrice}{ \begin{pmatrix} }{ \end{pmatrix} }
\newenvironment{pmatrice}{ \left( \begin{smallmatrix} }{ \end{smallmatrix}  \right) }

\newcommand{\gvettore}[1]{{ \begin{gmatrice} #1 \end{gmatrice} }}

\newcommand{\freccia}{{\rightarrow}}
\newcommand{\ifreccia}{{\hookrightarrow}}
\newcommand{\sfreccia}{{\rightarrow}}
\newcommand{\bfreccia}{{\hookrightarrow}}

\newcommand{\tende}{{\longrightarrow}}

\newcommand{\numero}{{\textsection}}

\newtheorem{tteo}{Theorem}[section]
\newtheorem{tprop}[tteo]{Proposition}
\newtheorem{tspe}{Hope}
\newtheorem{tlemma}[tteo]{Lemma}
\newtheorem{tcorol}[tteo]{Corollary}

\theoremstyle{definition}
\newtheorem{defin}[tteo]{Definition}

\newenvironment{teo}[1][]{   \begin{tteo}   {\bfseries \upshape #1} }{ \hspace{\stretch{1}} $\Box$ \end{tteo}   }
\newenvironment{prop}[1][]{  \begin{tprop}  {\bfseries \upshape #1} }{ \hspace{\stretch{1}} $\Box$ \end{tprop}  }

\newenvironment{lemma}[1][]{ \begin{tlemma} {\bfseries \upshape #1} }{ \hspace{\stretch{1}} $\Box$ \end{tlemma} }
\newenvironment{corol}[1][]{ \begin{tcorol} {\bfseries \upshape #1} }{ \hspace{\stretch{1}} $\Box$ \end{tcorol} }

\newcommand{\dimo}{  {\itshape Proof} \upshape :  }

\begin{document}

\sloppy

\maketitle

\begin{abstract}

In this paper we try to look at the compactification of Teichm\"uller spaces from a tropical viewpoint. Another paper working in this direction is \cite{FG}. Here we propose a completely different approach. We use amoebas and Maslov dequantization to construct and study the boundary of Teichm\"uller spaces, that can be seen as connected components of a real algebraic variety. 

We describe a general construction for the compactification of algebraic varieties, starting from their amoebas. This compactification is similar to the one described by Morgan e Shalen in \cite{MS1}, the difference is that they looked only to ``point at infinity'' to add, while we remove the points with some null coordinates and we replace them with some new points (see below for a definition). The boundary that we construct will be a closed subset of the sphere, such that the cone over this subset is a tropical variety.

When we apply this construction to the Teichm\"uller spaces we see that they can be mapped in a real algebraic hypersurface in such a way that the cone over the boundary is a subpolyhedron of a tropical hypersurface.   

We want to show how some properties of the boundary becomes straightforward if looked at from this point of view. For example there is a piecewise linear structure that appear naturally on the boundary, simply because the tropical varieties are polyhedrons. This structure is then shown to be equivalent to the one defined by Thurston in the 80's. 

Also we may see easily that every polynomial relation among trace functions on Teichm\"uller space may be turned automatically in a tropical relation among intersection forms over the boundary. This fact was already known (see \cite{Lu1} and \cite{Lu2}), but in this context it gain a theoretical justification.     

\end{abstract}
\ \\
\tableofcontents

\section{Amoebas}

In this section we recall the definition of the amoeba of an algebraic variety embedded in $\cappa^n$, and we state the fundamental facts about these objects that we will use in the following. Then we extend this definition to the case of an abstract variety. 

To define abstract varieties we fix a field $\cappa$ and a countable subfield $k \subset \cappa$. Let $V \subset \cappa^n$ and $W\subset \cappa^m$ be algebraic subvarieties defined over $k$. A polynomial map (with coefficients in $k$) $f:V\freccia W$ is said to be a polynomial isomorphism if it is bijective an if its inverse is again a polynomial map. By a $\cappa$-\nuovo{abstract affine variety} defined over $k$ we mean an equivalence class of such varieties up to polynomial isomorphisms.

This definition coincides with the usual one only if $\cappa$ is algebraically closed. Else this definition is more restrictive as we want a morphism to be polynomial. In this way every $\cappa$-abstract affine variety $V$ defined over $k$ has a well defined ring of coordinates $k[V]$. 

The hypothesis of $k$ being countable is not restrictive, as every variety $V \subset \cappa^n$ admits such a definition field. This hypothesis is useful in the sequel to make the ring of coordinate $k[V]$ countable.

\subsection{The tropical semifield and Maslov dequantization}

We define the \nuovo{tropical semifield} as the semifield $\erre^{trop} = (\erre,\oplus,\odot)$, where $a\odot b = a + b$, $a \oplus b = \max(a,b)$. It is called semifield as the $\odot$ operation is invertible ($0$ being the neutral element), while the $\oplus$ operation is not.

The Maslov dequantization is a continuous deformation of the semifield $\erre^{>0}$ to the semifield $\erre^{trop}$. Formally it is obtained in the following way. 

Let $h > 0$, and let $f_h = \exp(\frac{1}{h})$. The semifield operations on $\erre^{>0}$ induces, through the bijection $D_h:\erre^{>0} \ni x \bfreccia \log_{f_h}(x) \in \erre$, two operations on $\erre$: $\oplus_h$ and $\odot_h$. Explicitly $a \oplus_h b = D_h(D_h^{-1}(a) + D_h^{-1}(b)) = \log_{f_h}(f_h^a + f_h^b)$, and $a \odot_h b = D_h(D_h^{-1}(a) \cdot D_h^{-1}(b)) = a + b$. We will denote by $S_h$ the semifield $(\erre, \oplus_h, \odot_h)$, isomorphic to $\erre^{>0}$. 

If $h$ is small, the semifield $S_h$ is ``very similar'' to the tropical semifield $\erre^{trop}$, in the following sense. Let $a,b \in \erre$. Then $a \odot_h b = a \odot b$. Suppose $a \leq b$ so that $b= \max(a,b)$. Now $b \leq \log_f(f^a + f^b) \leq \log_f( 2 f^b ) = \log_f(2) + b$. Then $$a \oplus b \leq a \oplus_h b \leq a \oplus b + h \log(2)$$ 

We will write $S_0 = \erre^{trop}$, when needed.

One may try to use the maps $D_h$ to study real algebraic varieties. 
Let $Z \subset \erre^n$ be an algebraic variety. For the moment we look only at positive points, the set $Z_+ = Z \cap {(\erre^{>0})}^n$. We denote by $D_h^n$ the map $D_h$ applied to each factor of ${(\erre^{>0})}^n$, $D_h^n: {(\erre^{>0})}^n \bfreccia \erre^n$, and we denote by $Z_h$ the image $D_h^n(Z_+)\subset \erre^n$. 

Let $I \subset \erre[X_1 \dots X_n]$ be the ideal defining $Z$. An element $f \in I$ may be written in the form: 

$$f = \sum_{\omega \in \ze^n} a_\omega X^\omega $$

Where $X^\omega = \prod_{j =1 }^n X_j^{\omega_j}$, and the set  $A_f = \{\omega \in \ze^n \ |\ a_\omega \neq 0\}$ is finite. 

We want to separate the positive monomial of $f$ from the negative ones, so we define $A_f^+=\{ \omega \in A_f \ |\ a_\omega > 0\}$, and  $A_f^-=\{ \omega \in A_f \ |\ a_\omega < 0\}$. Now we split $f$ in its positive and negative part: 

$$f^+ = \sum_{\omega \in A_f^+} a_\omega  X^\omega \ \ \ ; \ \ \ f^- =\sum_{\omega \in A_f^-} (-a_\omega)  X^\omega $$

So $f = f^+ - f^-$. Then the set $Z_+$ may be written as:

$$ Z_+ = \{ x \in {(\erre^{>0})}^n \ |\  \forall f \in I : f^+(x)  = f^-(x)\}$$ 

So we have a set of equations, one for each polynomial in $I$, with positive coefficients. 

We are interested in the set $Z_h$, the image of $Z_+$ under the map $D_h^n$. For every polynomial $f \in I$ we take the transformation of its positive and negative parts through $D_h$. 

$$ f_h^+ = D_h \circ f^+ \circ {(D_h^n)}^{-1} = \begin{array}{c}\scriptstyle\ \\{\displaystyle\bigoplus}_h \\ {\scriptstyle\omega \in A_f^+}\end{array} D_h(a_\omega)\odot X^{\odot\omega}$$ 
$$f_h^- = D_h \circ f^- \circ {(D_h^n)}{-1} = \begin{array}{c}\scriptstyle\ \\{\displaystyle\bigoplus}_h \\ {\scriptstyle\omega \in A_f^-}\end{array} D_h(-a_\omega)\odot X^{\odot\omega}$$

As $D_h$ is a semifields isomorphism we have that 

\begin{equation}  \label{eqn:realvar}
Z_h = \{x \in \erre^n \ |\ \forall f \in I : f_h^+(x) = f_h^-(x)\}
\end{equation}

Then one may try to use Maslov dequantization to study the sets $V_h$.

\subsection{Absolute values}         \label{subsez:absolute values}

A similar technique also works for varieties over other fields. If $\cappa$ is a field we need a way for sending its multiplicative group in $\erre^{>0}$, then we may apply the maps $D_h$. The right way is an \nuovo{absolute value} function, i.e. a function   

$$ |\cdot|:\cappa \freccia \erre^{\geq 0} $$

satisfying: 1) $|x| = 0 \Leftrightarrow x = 0$;  2) $|xy| = |x||y|$; 3) $|x+y| \leq |x| + |y|$.

Let $\cappa$ be a field endowed with an absolute value $|\cdot|$. We may define the map $\Log$:

$$ \Log: {(\cappa^*)}^n \ni \gvettore{z_1 \\ \vdots \\ z_n} \freccia \gvettore{\log|z_1| \\ \vdots \\ \log|z_n|} \in \erre^n $$

This map is the composition of $D_1^n$ with the componentwise absolute value.

Two absolute values ${|\cdot|}_1, {|\cdot|}_2$ are said to be \nuovo{equivalent} if there exists a $\lambda\in \erre^{>0}$ such that ${|\cdot|}_1 = {|\cdot|}_2^\lambda$. Replacing $|\cdot|$ with an equivalent absolute value the $\Log$ map change by a scalar factor, i.e. as if we had used a map $D_h^n$ instead of $D_1^n$.

Let $V \subset \cappa^n$ an affine algebraic variety. We define the \nuovo{ameba} of $V$ as the set

$$\ameba(V) = \Log(V \cap  {(\cappa^*)}^n)$$

Let $I \subset \cappa[X_1 \dots X_n]$ be the ideal defining $V$. Again we choose an element $f = \sum a_\omega X^\omega \in I$. Every element $x \in V\cap  {(\cappa^*)}^n$ verifies:

$$\sum_{\omega \in A_f} a_\omega  x^\omega = 0$$

Using triangular inequality, this implies that

$$\forall w \in A_f : |a_\omega| |x|^\omega \leq \sum_{\omega \in A_f \setminus \{w\} } |a_\omega| |x|^\omega$$

where the symbol $|x| \in \erre^{>0}$ is a vector whose component are the absolute values of the components of $x$.

Hence every point $y \in \ameba(V) \subset \erre^n$ satisfy:

\begin{equation}   \label{eqn:complexvar}
\forall w \in A_f : D_1(|a_w|) \odot_1 y^{\odot_1 w} \leq \begin{array}{c}\ \\{\displaystyle\bigoplus}_1 \\ {\scriptstyle\omega \in A_f\setminus \{w\}}\end{array} D_1(|a_\omega|) \odot_1 y^{\odot_1\omega}
\end{equation}

One may ask which are the fields with an absolute value. First note that every field has the \nuovo{trivial} absolute value, i.e. an absolute value such that $\forall x \neq 0 : |x| = 1$. Excluding this case, the absolute values may be divided in archimedean ones and non archimedean ones.  

An absolute value is said to be \nuovo{archimedean} if  $\forall x, y \in \cappa : \exists n \in \enne : |nx| > |y|$. The best example is when $\cappa \subset \ci$ and $|\cdot|$ is the Euclidean norm. Actually every field with an archimedean absolute value may be embedded in $\ci$ in such a way that the norm is equivalent to the restriction of the Euclidean norm.  

The amoebas of varieties over archimedean fields are called \nuovo{archimedean amoebas}.

An absolute value $|\cdot|$ is non-archimedean if and only if it satisfies the \nuovo{ultrametric inequality}: 
$$ \forall x,y \in \cappa : |x+y|\leq \max(|x|,|y|) $$
Equivalently $|\cdot|$ is non archimedean if and only if the function
$$ v:\cappa^* \in x \freccia -\log(|x|) \in \erre $$
is a valuation.

On the other end, we may take any rank $1$ valuation on $\cappa$:
$$v:\cappa^* \freccia \Lambda$$
We have an immersion $\mu:\Lambda \ifreccia \erre$ well defined up to scalar multiplication. The map 
$$ |\cdot|:\cappa^* \ni z \freccia e^{-\mu(v(z))} \in \erre^+ $$
extended with $|0|=0$, is a non archimedean absolute value on $\cappa$, well defined up to equivalence.

When $|\cdot|$ is induced by a rank $1$ valuation, $v:\cappa^* \freccia \erre$, we can write the $\Log$ map as:
$$\Log: {(\cappa^*)}^n \ni \gvettore{z_1 \\ \vdots \\ z_n} \freccia \gvettore{-v(z_1) \\ \vdots \\ -v(z_n)} \in \erre^n$$

And the amoeba becomes
$$\ameba(V) = \left\{\gvettore{-v(z_1) \\ \vdots \\ -v(z_n)} \ |\ \gvettore{z_1 \\ \vdots \\ z_n} \in V \cap  {(\cappa^*)}^n\right\} $$

An amoeba of this kind is said to be a \nuovo{non-archimedean amoeba}.

As the ultrametric inequality is much stronger than the triangular one, the relation \ref{eqn:complexvar} becomes:

\begin{equation}  \label{eqn:nonarchvar}
\forall w \in A_f : (-v(a_w)) \odot_1 y^{\odot_1 w} \leq \max_{\scriptstyle\omega \in A_f\setminus \{w\}} (-v(a_\omega)) \odot_1 y^{\odot_1\omega} 
\end{equation}

\subsection{Tropical varieties}

By a tropical polynomial in $n$ variables we mean an expression of the form 

$$P(X) = \bigoplus_{\omega \in A_P} a_\omega\odot X^{\odot \omega}$$

Here $A_P \subset \ze^{\oplus J}$ is a finite set, $a_\omega \in \erre$ and $X^{\odot \omega} = \bigodot_{j = 1}^n X_j^{\odot \omega_j}$. Such a polynomial defines a function: 

$$P:\erre^n\ni x \freccia P(x) =  \max_{\omega \in A_P}(a_\omega + <x,\omega>) \in \erre$$

This is a convex piecewise linear function. 

Tropical polynomials may be used to define tropical hypersurfaces. We see two different way for doing this.

If $P$ and $Q$ are tropical polynomials with  $A_P \cap A_Q = \emptyset$, we define the \nuovo{real tropical hypersurface} defined by $P$ and $Q$ as 

\begin{equation}
T_\erre(P,Q) = \{x\in \erre^n \ |\ P(x)=Q(x) \}
\end{equation}

This definition is motivated by the relation \ref{eqn:realvar}. Note that the hypothesis $A_P \cap A_Q = \emptyset$ is perfectly natural in this context, as we think $P$ and $Q$ to be the positive and negative part of one polynomial, and it is necessary to guarantee that the set $T_\erre(P,Q)$ has empty interior.

An hypersurface $T_\erre(P,Q)$ is a union of (possibly noncompact) polyhedrons, and in many cases it is a topological manifold. 

This kind of varieties was used in \cite{Vi} to give a tropical account of the patchworking theorem, and will be used in section \ref{subsez:real comp} to describe the boundary we construct for a real algebraic variety.

The most common way for defining hypersurfaces is inspired by relation \ref{eqn:complexvar}. If $P$ is a tropical polynomial, we denote by $T(P)$ the \nuovo{(complex) tropical hypersurface} associated with $P$, defined to be the set points in $\erre^n$ satisfying the following relation:

\begin{equation}
\forall w \in A_P : a_w\odot X^{\odot w} \leq \bigoplus_{\omega \in A_P \setminus \{w\} } a_\omega\odot X^{\odot \omega}
\end{equation}

This may be expressed more concisely as the set points in $\erre^n$ in which at least two of the monomials of $P$ achieve maximum, or equivalently as the set of points in $\erre^n$ in which the function $P$ is not locally linear. Tropical hypersurfaces are closely related with non-archimedean amoebas for hypersurfaces:

\begin{teo}        \label{teo:tropvar}
Let $\cappa$ be a algebraically closed field endowed with a valuation $v:\cappa \freccia \erre\cup\{\infty\}$, let $f = \sum_{\omega \in A_f} a_\omega X^\omega \in \cappa[X_1 \dots X_n]$, and let $V\subset \cappa^n$ the zero locus of $f$. Then the closure of the amoeba $\ameba(V)$ coincides with the tropical hypersurface defined by the tropical polynomial $f^\tau = \bigoplus_{\omega \in A} (-v(a_\omega))\odot X^{\odot \omega}$.

\dimo See \cite[thm.~2.1.1]{EKL}. Actually we have already proved that $\ameba(V) \subset T(P)$, as it follows easily from relation \ref{eqn:nonarchvar}, if one remember that $\odot_1 = \odot$, and $\max = \oplus$. The reversed inequality follows from \cite[lemma~2.1.5]{EKL}.
\end{teo}

This result may be used as a guide in the definition of \nuovo{(complex) tropical variety} of lower dimension. We choose a field $\cappa$ with a valuation $v$, we take an ideal $I \subset \cappa[X_1 \dots X_n]$, and we define the tropical variety of the ideal $I$ as:

\begin{equation}
T(I) = \bigcap_{f \in I} T(f^\tau)
\end{equation}
  
Were $f^\tau$ is the tropical polynomial defined in theorem \ref{teo:tropvar}. With this definition we have that every tropical variety is an intersection of tropical hypersurfaces, but is not true that an arbitrary intersection of tropical hypersurfaces is a tropical variety.

\begin{teo}
Let $\cappa$ be an algebraically closed field endowed with a valuation. Let $I \subset \cappa[X_1 \dots X_n]$ be an ideal, and let $V \subset \cappa^n$ be its zero locus. Then the closure of the amoeba of $V$ is equal to $T(I)$.

\dimo See \cite[thm. 2.1]{SS}. 
\end{teo}

Finally there is a third description for tropical varieties, using ring valuations. 
\begin{defin}
Let $k$ be a field, $A \supset k$ a $k$-algebra. By a real valued \nuovo{valuation} of $A$ over $k$ we mean a map $v:A \freccia \erre \cup \{+\infty\}$ such that: 
\begin{enumerate}
	\item $v(xy) = v(x) + v(y)$.
	\item $v(x + y) \geq \min(v(x),v(y))$.
	\item $v(0) = \infty$, $v(k^*) = 0$.
\end{enumerate}

A valuation is said to be \nuovo{trivial} if $v(A) \subset \{0,+\infty\}$.
\end{defin}

Let $\cappa$ be an algebraically closed non-archimedean field, with valuation $v$, and $X \subset \cappa^n$ be an irreducible variety. Let $A=\cappa[X\cap {(\cappa^*)}^n]$ and $S_\erre$ be the set of real valuation of $A$ extending $v$. The Bieri-Groves set of $X$ is defined as the image of $S_\erre$ by the map 
$$z:S_\erre \ni w \freccia (-w(X_1) \dots -w(X_n)) \in \erre^n$$

\begin{teo} \label{teo:bierigroves}
$z(S_\erre)=\ameba(X)$.

\dimo See \cite[thm.~2.1.1]{EKL}  
\end{teo}

There is a property of valuations that we will need in the following. The set $v^{-1}(+\infty)$ is a prime ideal of $A$. Hence all invertible elements are real valued. If $I \subset A$ is an ideal and $v$ is a valuation of the quotient $A/I$, then $v$ may be lifted to a valuation $w$ of $A$, by composition with the quotient map. Then $w(I) = \{+\infty\}$. The reciprocal is also true:

\begin{lemma}      \label{lemma:quot val}
Let $w$ be a valuation of a ring $A$. Suppose $I$ is an ideal such that $w(I) = \{+\infty\}$. Then there exists a valuation $v$ of $A/I$ such that $w$ is equal to the composition of $v$ with the quotient map.

\dimo We only need to say that if $f,g \in A$ and $h=f-g \in I$, then $w(f)=w(g)$. If both have infinite valuation this is obvious, else we suppose that $f$ has finite valuation and we know that $f + h = g$. As $w(f)<w(h)$ we have that $w(g)=w(f)$. So we may define a valuation $v$ of $A/I$, taking the valuation of a counterimage. 
\end{lemma}

\subsection{Amoebas for abstract varieties}   \label{subsez:amebevarastratte}

Now we want to extend the definition of amoeba to the case of abstract varieties. 

Let $\cappa$ be a field, $V$ an abstract affine variety defined over the countable subfield $k$. If we choose any immersion $V \subset \cappa^n$, the family of the coordinate functions $\{X_1 \dots X_n\}$ generates the ring $k[V]$ as a $k$-algebra.
 
Reciprocally, let $\famil = \{f_1 \dots f_n\}$, $f_i \in k[V]$, be a fixed family whose elements generate the ring $k[V]$ as a $k$-algebra. The family $\famil$ induces a map $\overline{\famil}$:

$$ \overline{\famil}:V \ni x \freccia \gvettore{f_1(x) \\ \vdots \\ f_n(x)} \in \cappa^n $$

\begin{prop}          \label{prop:funzioni generatrici}
The map $\overline{\famil}$ is injective, its image $\overline{\famil}(V)$ is an algebraic subvariety and $\overline{\famil}$ is a polynomial isomorphism between $V$ and $\overline{\famil}(V)$.

\dimo Injectivity follows easily as $\famil$ generate $k[V]$. The remaining may be done as in the proof of \cite[I\numero 3, prop.1]{M2}
\end{prop}

If we identify $V$ with $\overline{\famil}(V)$, the function $f_i$ becomes the coordinate function $X_i$. This shows that choosing an immersion in $\cappa^n$ is equivalent to choosing a family $\famil = (f_1 \dots f_n)$, $f_i \in k[V]$, whose elements generate the ring $k[V]$ as a $k$-algebra. 

This construction may be generalized: we may take a countable family $\famil = {(f_j)}_{j \in J}$, $f_j \in k[V]$, with $\card(J) \leq \aleph_0$ such that its elements generate $k[V]$ as a $k$-algebra. Such a family will be called a \nuovo{generating family}.

A generating family $\famil = {(f_j)}_{j \in J}$ defines an injective and proper map:

$$\overline{\famil}: V \ni x \freccia {(f_j(x))}_{j \in J} \in \cappa^n$$

Now suppose that the field $\cappa$ is endowed with an absolute value. Let $\famil$ be a fixed generating family. We define $V' = \{x\in V \ |\ \forall j \in J : f_j(x) \neq 0 \}$. Note that $V'$ may be empty (this may happen trivially if $0 \in \famil$, but it may also happen sometimes when $V$ is reducible). The set $V \setminus V'$ is a countable union of Zariski closed set of $V$. We define the $\Log$ map:

$$ \Log:V' \ni x \freccia {(\log|f_j(x)|)}_{j \in J} \in \erre^J$$

If $\famil$ is finite $\Log = \Log \circ \overline{\famil}|_{V'}$, hence $\ameba(\overline{\famil}(V)) = \Log(V')$.

In the general case we may define the amoeba of $V$ as 

$$\ameba(V) = \Log(V') \subset \erre^J$$

In this definition $V'$ plays the role that in the previous section was played by $V \cap {(\cappa^*)}^n$. 

The ideal of a variety is usually defined only for subvarieties of $\cappa^n$. Actually we can define it also for an abstract affine variety with a fixed generating family. 

We need the free $k$-algebra generated by the symbols ${(X_j)}_{j\in J}$, that will be written as $k[(X_j)]$. If $J$ is finite it is a usual polynomial ring, else it is the ring of polynomials in an infinite number of variables. There exists a unique $k$-algebras homomorphism $\varphi:k[(X_j)] \sfreccia k[V]$ such that $\varphi(X_j) = f_j$. This homomorphism is surjective, as the elements $f_j$ generate $k[V]$. We define the \nuovo{ideal} of $V$ with reference to the family $\famil$: $I = \ideal(V) = \ker \varphi$. This implies $k[(X_j)] / I \simeq k[V]$, hence $I$ is a prime ideal if and only if $V$ is irreducible. 

We want to define the \nuovo{ring of coordinates} of $V'$. If $V$ is irreducible we may simply define it as $k[V'] = k[(f_j), (f_j^{-1})] \subset k(V)$. If $V$ is reducible we denote by $S$ the multiplicative set generated from $\famil$. Then we define $k[V'] = S^{-1}k[V]$, the localization with reference to $S$. Note that $k[V']=0$ if and only if $0 \in S$ if and only if $V'=\emptyset$. 

Now we may define the \nuovo{ideal} of $V'$. We need the $k$-algebra of Laurent polynomials, $k[(X_j), (X_j^{-1})] = k[{(X_j)}_{j \in J}, {(Y_j)}_{j \in J}] / ({\{X_j Y_j - 1\}}_{j \in J})$. There exists a unique homomorphism $\psi:k[(X_j),(X_j^{-1})] \sfreccia k[V']$ such that $\psi(X_j) = f_j$. We define the ideal of $V'$ as before: $I'=\ideal(V')=\ker \psi$. If $\famil$ is finite $V'$ is an affine variety, else it may not, but, in a formal sense, it has a defining ideal.

In this more general context we can give the definition of tropical polynomials and of tropical varieties. Let $J$ be a (finite or) countable set. By a tropical polynomial in the variables ${\{X_j\}}_{j\in J}$ we mean an expression of the form 

$$P(X) = \bigoplus_{\omega \in A_P} a_\omega\odot X^{\odot \omega}$$

Here $A_P \subset \ze^{\oplus J}$ is a finite set,  $a_\omega \in \erre$ and $X^{\odot \omega} = \bigodot_{j \in J} X_j^{\odot \omega_j}$, a finite product as all $\omega_j$ but a finite number are $0$. Such a polynomial defines a convex piecewise linear function (for the linear structure of $\erre^J$, see below): 

$$P:\erre^J\ni x \freccia P(x) =  \max_{\omega \in A_P}(a_\omega + <x,\omega>) \in \erre$$

Now one can easily give definition for real and complex tropical hypersurfaces in $\erre^J$, $T_\erre(P,Q)$, $T(P)$, similar to the definitions given in previous subsection. 

\section{A Morgan-and-Shalen-like compactification}

In this section we construct a compactification for complex and real varieties that is similar to the compactification constructed by Morgan and Shalen in $\cite{MS1}$. The compactification presented here is more suitable to be described by tropical objects, as it is constructed starting from the $\Log$ map. 

\subsection{Infinite-dimensional spaces}

Let $J$ be a (finite or) countable set. We need some properties of the space $\displaystyle \erre^J = \prod_{j \in J} \erre$. If endowed with the product topology it is a second countable, metrizable Hausdorff space. If $J$ is infinite it isn't locally compact, but the closed subsets of the form $\displaystyle \prod_{j \in J} [\alpha_j, \beta_j]$ are compact. 

\begin{prop}
Let $V$ be a real or complex variety, and let $\famil= {\{f_j\}}_{j \in J}$ be a generating family for $V$. Then the amoeba $\ameba(V)$ is a closed subset of $\erre^J$.

\dimo As $\erre^J$ is second countable it is enough to show that if a sequence $(y_i) \subset \ameba(V)$ converges to a point $y \in \erre^J$, then $y \in \ameba(V)$. Choose a sequence $(x_i) \subset V'$ such that $\forall i: \Log(x_i) = y_i$. The sequence $(x_i)$ is contained in a compact subset of $V$, else we could find a subsequence $x_{k_i}$ and an element $f_j\in\famil$ such that $|f_j(x_{k_i})| \tende \infty$ but this implies that $\Log(x_i)$ may not converge. Hence there is a subsequence $(x_{h_i})$ converging to some point $x \in V$. We know that $x \in V'$, else we could find an element $f_j\in\famil$ such that $f_j(x)=0$, and again $\Log(x_i)$ may not converge. Hence we have $y_{k_i} \tende y = \Log(x) \in \ameba(V)$.
\end{prop}

The space $\erre^J$ has a natural structure of real vector space. The elements $e_i = {(\delta_{ij})}_{j \in J}$ are independent and they generate a subspace that will be denoted by $\erre^{\oplus J} = \displaystyle \bigoplus_{j \in J} \erre$. If $J$ is finite this subspace coincides with the whole $\erre^J$, else it is smaller. We have a natural pairing:
 
$$ \erre^{J} \times \erre^{\oplus J} \in (\xi,x) \freccia <\xi,x> = \sum_{j\in J} \xi_j x_j \in \erre $$

It is well defined as all elements of the sum but a finite number are $0$. By \nuovo{coordinate functions} on $\erre^J$ we mean the functionals $e^i:\erre^J \ni \xi \freccia <\xi,e_i> \in \erre$.

We denote by $\sim$ the following equivalence relation on $\erre^J \setminus \{0\}$: $x \sim y$ $\Leftrightarrow$ $\exists \alpha \in \erre^{>0} : x = {(t_j)}_{j \in J}$ and $y = {(\alpha t_j)}_{j \in J}$. The quotient $\sferic = (\erre^J \setminus \{0\}) / \sim$ is a standard sphere with dimension $|J|-1$ when $J$ is finite, else it is a sort of infinite-dimensional sphere. In either case it is a metrizable Hausdorff space. We will denote the projection by $\pi:\erre^J \setminus \{0\} \sfreccia \sferic$. Points of $\sferic$ will be denoted by homogeneous coordinates: ${[t_j]}_{j \in J}=\pi({(t_j)}_{j \in J})$.

Suppose $C \subset \displaystyle \prod_{j \in J} [\alpha_j, \beta_j]$ is a closed set not containing $0$, it is a compact subset of $\erre^J \setminus \{0\}$, hence $\pi(C)$ is compact. However, if $J$ is infinite, there isn't any such set $C$ satisfying $\pi(C) = \sferic$, as $\sferic$ is not compact.

Now we glue together $\erre^J$ and $\sferic$ obtaining $D^J$, a \nuovo{closure at infinity} of $\erre^J$. If $J$ is finite $D^J$ is simply the compactification of $\erre^J$ with a sphere at infinity, while in the general case it is not compact.  

As a set $D^J$ is $\erre^J \cup \sferic$. Then we define a basis for the topology. Let $U$ be an open subset of $\sferic$. A subset $V$ is said to be an \nuovo{open conic subset} with base $U$ if it may be written as $V = U \cup (\pi^{-1}(U) \cap H_1 \cap \dots \cap H_n)$, where $H_i$ are coordinate semispaces, i.e. sets of the form $\{x \in \erre^J \ |\ e^h(x) > c\}$ o $\{x \in \erre^J \ |\ e^h(x) < c\}$, with $h\in J$, $c \in \erre$. Let $B', B''$ be basis for the topology of $\erre^J$ and $\sferic$ respectively. A basis for the topology of $D^J$ is given by:

$$ B' \cup \{  V  \ |\ V \mbox{ is an open conic subset with base } U \in B'' \} $$

With reference to this topology $\erre^J$ is a dense open subset, and $\sferic$ is a closed subset with empty interior.

\begin{prop}  \label{prop:sfera}
A sequence $(x_n) \subset \erre^J$ converges to a point $x\in\sferic$ if and only if:
\begin{enumerate}
	\item $\exists j \in J : e^j(x_n) \tende +\infty$ o $e^j(x_n) \tende -\infty$ in $\erre$.
	\item For every subsequence $(x_{n_k}) \subset \erre^J \setminus \{0\}$ we have $\pi(x_{n_k}) \tende x$ in $\sferic$.
\end{enumerate}

\dimo It follows from the definition of the topology.
\end{prop}

When $J$ is finite, $|J|=n$, we could have constructed $D^J$ in a more concrete way:

$$ \phi:\erre^n \ni x \ifreccia \frac{x}{\sqrt{1+{\|x\|}^2}} \in D^n$$

The image of this map is the interior of $D^n$, and it sends rays of $\erre^n$ in rays of $D^n$ and spheres of $\erre^n$ centered in the origin in spheres of $D^n$. The map $\phi$ may be extended to an homomorphism from $D^J$ to $D^n$.

\subsection{Compactifications determined by a map}        \label{sez:compatt}

Let $X$ be a locally compact Hausdorff topological space. By a \nuovo{compactification} of $X$ we mean a map $\eta:X \freccia Y$, where $Y$ is a compact Hausdorff space, and $\eta$ is a homeomorphism of $X$ into a dense open subset of $Y$. We denote the Alexandrov compactification of $X$ by $\Hat{X} = X \cup \{\infty\}$.

Let $P$ be an Hausdorff space and $t:X \freccia P$ be a continuous map whose image is relatively compact in $P$. We recall a construction given in \cite{MS1} for a compactification of $X$ that will be called compactification \nuovo{determined} by $t$.  

Let $i$ be the map:

$$i:X \ni x \freccia (x,t(x)) \in \Hat{X} \times P$$

The codomain of the compactification is $Y = \overline{i(X)} \subset \Hat{X} \times P$, and the map is $\eta = i|_X^Y$. The space $Y$ is a closed subset of $\Hat{X} \times \overline{t(X)}$, hence it is compact.

By $p_1:\Hat{X} \times P \freccia \Hat{X}$ and $p_2:\Hat{X} \times P \freccia P$ we denote the canonical projection on the first and the second factor, respectively.

The set $B = {(p_1|_Y)}^{-1}(\infty) \subset Y$ is a closed subset with empty interior, and its complement $Y \setminus B$ is equal to $\eta(X)$, hence $\eta(X)$ is a dense open subset of $Y$. This shows that $\eta$ is a compactification of $X$. The set $B$ is said to be the set of \nuovo{ideal points} of the compactification.

The projection $p_2|_B$ is a homeomorphism with its image, (as $B\subset\{\infty\}\times P$), hence we may identify $B$ with a subset of $P$. We may also identify $X$ with $\eta(X)$, and so we may think of $Y$ as $X \cup B$. A sequence $\{x_n\} \subset X$ converges to a point $b \in B$ if and only if $x_n \tende \infty$ in $\Hat{X}$ and $t(x_n) \tende b$ in $B$.

If the spaces $X$, $\Hat{X}$ e $P$ are first countable this fact characterizes the topology of $Y$, and we may describe $B$ as: $B = \{ b \in P \ |\ \exists \{x_n\} \subset X: x_n \tende \infty$ in $\Hat{X}$ e $t(x_n) \tende b$ in $P \}$.

\subsection{Compactification of complex varieties}

Let $V$ be a complex abstract affine variety defined over the countable field $k \subset \ci$, embedded with the classical topology, and let $\effe = k(V)$. We fix a generating family $\famil \subset k[V]$, $\famil = {(f_j)}_{j \in J}$ and we denote by $V'$ the subset $\{x \in V \ |\ \forall j \in J : f_j(x) \neq 0 \}$. 

We denote by $\theta:V' \freccia D^J$ the map $\Log$ when thought as a map taking values in $D^J$, the closure at infinity of $\erre^J$

Let $V_1 \dots V_n$ be the irreducible components of $V$. For every component $V_i$ there is a prime ideal $P_i \subset k[V]$ such that $k[V_i] = k[V] / P_i$. We define the set $\famil_i$ as set of the images in $k[V_i]$ of the elements of $\famil$. The set $\famil_i$ is a generating family for $V_i$. With reference to this family we define the set $V_i'$ and the ring $k[V_i']$. The $\Log$ map for $V_i'$ is exactly the restriction of the $\Log$ map for $V'$.

To take the compactification defined by $\theta$ we need that the closure of the image $\overline{\theta(V')} \subset D^J$ is compact. This is obvious when $J$ is finite.

\begin{prop}
The set $\overline{\theta(V')} \subset D^J$ is compact.

\dimo Without loss of generality we may add the hypothesis that $V$ is irreducible, as compactness is perserved by finite unions. As $D^J$ is second countable we only need to prove that every sequence $(y_n) \subset \theta(V')$ has a subsequence that converges to some point $y \in D^J$. Take a sequence $(x_n) \subset V'$ such that $\forall i:\theta(x_i)=y_i$. If $x_n$ is infinitely often in a compact subset of $V'$ we can extract a subsequence $(x_{k_n})$ converging to some point $x \in V'$, so that the subsequence $(y_{k_n})$ converges to $y=\theta(x) \in \theta(V')$. Else every compact subset of $V'$ contains only a finite number of elements of $x_n$. We conclude using some notation and proposition that will be introduced in subsection \ref{subsez:assocval}. From $x_n$ we may extract a quasi-valuating subsequence $(x_{h_n})$. By proposition \ref{prop:compatibile} there exists a valuation $v$ compatible with $(x_{h_n})$. By lemma \ref{lemma:limite} we may conclude that the subsequence $y_{h_n} = \theta(x_{h_n})$ has limit in $D^J$. 
\end{prop}

We denote by ${V'}^{comp}=V' \cup B(V')$ the compactification of $V'$ determined by $\theta$. The set $B(V')$ may be seen as a subset of $D^J$:

$$B(V') = \{ b \in D^J \ |\ \exists \{x_n\} \subset V': x_n \tende \infty \mbox{ in } \Hat{V'} \mbox{ e } \theta(x_n) \tende b \mbox{ in } D^J \}$$

$x_n \tende \infty$ in $\Hat{V'} \Rightarrow \|\Log(x_n)\| \tende \infty \Rightarrow \|\theta(x_n)\| \tende \infty$. Hence $B(V') \subset \sferic$.

Recall that $\pi:\erre^J \setminus \{0\} \sfreccia \sferic$ is the projection in the quotient, we get:

$$ B(V') = \{ b \in \sferic \ |\ \exists \{x_n\} \subset V': x_n \tende \infty \mbox{ in } \Hat{V'} \mbox{ e } \pi(\Log(x_n)) \tende b \mbox{ in } \sferic \} $$
 
If $0 \not\in \Log(V')$ we could have defined $\theta = \pi \circ \Log:V' \freccia \sferic$, and we had gotten the some compactification with the same set of ideal points $B(V')$.

When $J$ is finite we may use the ``concrete'' description of $D^J$ as $D^n$:

$$ \phi:\erre^n \ni x \ifreccia \frac{x}{\sqrt{1+{\|x\|}^2}} \in D^n$$

The compactification is defined by the map $\theta=\phi \circ \Log$, exactly the map used by \cite{Be2}.

\subsection{Compactification of real varieties}   \label{subsez:varreali}

Let $Z \subset \erre^n$ be a real variety defined over the countable subfield $k \subset \erre$, and let $Z_\ci$ be its complexification, i.e. the smaller complex variety in $\ci^n$ containing $Z$. We may choose $k$ such that $Z_\ci$ is defined over $k$. We have that $k[Z_\ci]=k[Z]$ and $Z_\ci \cap \erre^n = Z$. 

From now on we fix a generating family $\famil = {(f_j)}_{j\in J} \subset k[Z]$ for $Z$, and we have $Z' = \{x \in Z \ |\ \forall j \in J : f_j(x) \neq 0 \}$. The family $\famil$ is also a generating family for $Z_\ci$, so it defines a compactification ${Z_\ci'}^{comp} = Z_\ci' \cup B(Z_\ci')$. 

We define the compactification $Z'^{comp}$ of $Z'$ to be the closure of $Z'$ in ${Z_\ci'}^{comp}$. Hence we have $B(Z')=\overline{Z'}\cap B(Z_\ci')$,  ${Z'}^{comp} = \overline{Z'} = Z' \cup B(Z')$. The set $B(Z')$ may be described explicitly as  $\{ b \in \sferic \ |\ \exists \{x_n\} \subset Z': x_n \tende \infty \mbox{ in } \Hat{Z'} \mbox{ and } \pi(\Log(x_n)) \tende b \mbox{ in } \sferic \}$.

Actually we need something more, we want to define a compactification for a connected component of a real variety. So let $C$ be a connected component of $Z'$, we take as compactification its closure in $Z'^{comp}$, so that $B(C) = \overline{C} \cap B(Z')$ and $C^{comp} = \overline{C} = C\cup B(C)$. Explicitly $B(C) = \{ b \in \sferic \ |\ \exists \{x_n\} \subset C: x_n \tende \infty \mbox{ in } \Hat{C} \mbox{ and } \pi(\Log(x_n)) \tende b \mbox{ in } \sferic \}$. 

So for a real variety we have a boundary that is contained in the boundary of its complexification. 

\subsection{Extension to non-generating families} \label{subsez:non generating}

We may easily construct compactifications also from some finite non-generating family $\famil \subset k[V]$.

Let $V$ be a complex affine variety defined over the countable subfield $k\subset \ci$. We choose a finite family $\famil = \{f_1 \dots f_n\}$, and we define the map

$$ \overline{\famil}: V \ni x \freccia \gvettore{f_1(x)\\ \vdots \\ f_n(x)} \in \ci^n $$

Let $V_1 \dots V_n$ the irreducible components of $V$. We denote by $W_i$ the closure (with reference to the Zariski topology) of the image of $V_i$ through the map $\overline{\famil}$. The closed sets $W_i$ are affine irreducible varieties (see \cite[I\numero 8, prop.~1]{M2}), and are defined over $k$. The map $\overline{\famil}$, thought as a map from $V_i$ to $W_i$ is a dominating map, hence its image contain a non empty, Zariski-open, subset $U_i \subset W_i$ (see \cite[VIII\numero 3, thm.~3]{M2}). 

Every $U_i$ is dense (with reference to the classical topology) in $W_i$ (see \cite[I\numero 10, thm.~1]{M2}), hence if we denote by $W$ the closure (in the classical topology) of the image of $\overline{\famil}$ we have that $W=W_1 \cup \dots \cup W_n$. Hence $W$ is a complex affine variety defined over $k$, and $\overline{\famil}$ may be seen as a dominating map from $V$ to $W$, whose image is dense in the classical topology. So $\overline{\famil}$ induces an identification of $k[W]$ with a subring of $k[V]$, the $k$-algebra generated by the functions $f_1 \dots f_n$. 

We choose the family $\{X_1 \dots X_n \}$ as a generating family for $W$. We may construct the compactification of $W'$ with reference to this family: $W'^{comp}=W' \cup B(W')$. The map $\overline{\famil}$ may be seen as a map from $V' = \{x\in V \ |\ \forall i \in \{1\dots n\} : f_i(x) \neq 0 \}$ to $W'$, and also as a map $\overline{\famil}:V' \freccia W'^{comp}$. This way the image of $\overline{\famil}$ is relatively compact, so we may consider the compactification of $V'$ determined by $\overline{\famil}$, and we denote it again by $V'^{comp} = V' \cup B(V')$. 

If we want to study the set of ideal points $B(V')$ through the methods presented above, we need some additional hypothesis on $\famil$.  

If we are interested to the whole $B(V')$ the hypothesis we need is that $\overline{\famil}:V'\freccia W'$ is a proper map. With this hypothesis the set $B(V') \subset W'^{comp}$ may be described as $B(V') = \{ b \in W'^{comp} \ |\ \exists \{x_n\} \subset V': x_n \tende \infty$ in $\Hat{V'}$ and $\overline{\famil}(x_n) \tende b$ in $W'^{comp} \}$. By properness if $\{x_n\}\subset V'$ is a sequence the property $x_n \tende \infty$ in $V'$ is equivalent to $\overline{\famil}(x_n) \tende \infty$ in $W'$. Hence $B(V') = B(W')$.

In the following we are interested in the compactification of a connected components of a real variety. Suppose know that $V$ is a complex variety defined over the countable field $k\subset \erre$, and that $C$ is a connected components of the real part $V_\erre'$. We may define a compactification $C^{comp}=C \cup B(C)$ as the closure of $C$ in $V'^{comp}$, $B(C)$ being $C^{comp} \cap B(V')$. To study $B(C)$ through the methods presented above we need a weaker hypothesis, simply we need $\overline{\famil}:C \freccia W'$ to be a proper map. With this hypothesis we have that $\overline{\famil}(C)$ is contained in a connected components $D$ of $W_\erre'$, and that $B(C) \subset B(D) \subset B(W')$.

\section{Tropical description}

In this section we want to see how Maslov dequantization may be useful to describe the compactification defined above.

\subsection{Compactification in the real case}      \label{subsez:real comp}

If $V \subset \erre^n$ and $\famil = \{X_1 \dots X_n\}$, we may decompose $\erre^n$ in its $2^n$ ``orthants'', i.e. if we fix an $s \in {\{-1,1\}}^n$, we define $\erre^n_s = \{x \ |\ \forall i: x_i s_i > 0 \}$. Then we consider the induced decomposition on $V$: $V_s = V \cap \erre^n_s$.
In the general case, when $\famil={\{f_j\}}_{j \in J}$ is a generic generating family, this corresponds to decompose $V$ in pieces, each of which indexed by an element $s \in {\{-1,1\}}^J$: $V_s = \{x\in Z \ |\ \forall j: f_j(x)s_j > 0\}$. In the following we fix such an $s$.

For every $h>0$ we define the map $D_h^J$:

$$D_h^J:V' \ni x \freccia {( D_h(|f_j(x)|) )}_{j \in J} \in \erre^J$$

The image $D_h^J(V_s)$ will be denoted by $V_{h,s}$. The map $D_1^J$ is the $\Log$ map. Then we use these maps to construct a dequantifying deformation of $V_s$:

$$D:V_s \times (0,+\infty) \ni (x,h) \freccia  (D_h^J(x),h)\in\erre^J \times (0,+\infty)$$

This map is a homeomorphism of the product $V_s \times (0,+\infty)$ in $W = \bigcup_{h > 0} V_h \times \{h\} \subset \erre^n \times (0,+\infty)$. So the image set $W_s$ may be viewed as a deformation of $V_s$, the deformation induced by Maslov dequantization. The map is proper so $W_s$ is closed. Let $\overline{W_s}$ be the closure of this set in $\erre^J \times [0,+\infty)$. As $W_s$ is a closed subset of $\erre^J \times (0, +\infty)$ we have that $\overline{W_s}$ is simply the union of $W_s$ with some subset of $\erre^J \times \{0\}$, that will be called $V_{0,s}$.  

\begin{prop}       \label{prop:cone}
The set $V_{0,s}$ is a cone.

\dimo This is easy as soon as one notes that $\forall x \in \erre^{>0}: D_h(x) = h D_1(x)$. If $x \in V_0$ there exists a sequence $(x_n, h_n) \subset W$ converging to $x$. 

For every positive real number $\lambda$ the sequence 
$$y_n = (D_{\lambda h_n}^J( {(D_{h_n}^J)}^{-1}(x_n) ), \lambda h_n) \in V_{\lambda h_n} \times \{\lambda h_n\}$$ 
is in $W$, and as $D_{\lambda h_n}^J( {(D_{h_n}^J)}^{-1}(x_n) = \lambda x_n$ it tends to $\lambda x$.
\end{prop}

Now we relate the set $V_{0,s}$ with the boundary constructed in the previous section. We denote by ${V_s}^{comp}$ the closure of $V_s$ in $V'^{comp}$, and we define $B(V_s)=B(V') \cap {V_s}^{comp}$.

\begin{prop}  \label{prop:boundary}
The cone $V_{0,s}$ is equal to $CB(V_s)$, the cone over the boundary.

\dimo $V_{0,s} \subset CB(V_s)$: If $x\in V_{0,s}$, $x \neq 0$, there exists a sequence $(x_n, h_n) \subset W$ converging to $x$. We pull it back on $V_s$ in this way: $y_n = {(D_{h_n | V_s}^J)}^{-1}(x_n) \in V_s$. Now $\Log(y_n)=D_1^J(y_n) = \frac{1}{h_n}x_n$ is a sequence tending to infinity in $V_s$. The sequence of projections in the sphere, $\pi(\Log(y_n)) = \pi(x_n)$, also converges to $\pi(x)$, so $\pi(x) \in B(V_s)$.

$CB(V^+) \subset V_0$: If $x \in CB(V^+)$, $x \neq 0$, there exists a sequence $(y_n) \subset V_s$ converging to $\pi(x)$. As $x \neq 0$ we can find an index $j$ such that the $j$-th coordinate of $x$ is equal to $\xi \neq 0$. Let $z_n$ be the $j$-th coordinate of $D_1^J(y_n)$. We may suppose $z_n \neq 0$. Let $h_n = \frac{\xi}{z_n}$, a sequence converging to $0$, as $z_n$ tends to infinity. The sequence $(D_{h_n}^J(y_n), h_n) \in W_s$ converges to $x$, so $x \in V_{0,s}$.  
\end{prop}

Now we try to give a description of $CB(V_s)$, using this information. First we recall a fact about tropical polynomials.

\begin{prop}
Let 
$$P = \bigoplus_{\omega \in A_P} a_\omega \odot X^{\odot\omega}$$ 
be a tropical polynomial, and let 
$$P_h = \begin{array}{c}{\scriptstyle \ } \\ {\displaystyle\bigoplus}_h \\ {\scriptstyle\omega \in A_P}\end{array} a_\omega \odot_h X^{\odot_h\omega}$$ 
be the corresponding family of polynomials in $S_h$. Then the family $P_h$ converges uniformly to $P$.

\dimo This follows from the inequality $a\oplus b \leq a \oplus_h b \leq a \oplus b + \log(2)$ that implies 
$$\forall x \in \erre : P(x) \leq P_h(x) \leq P(x) + h \log(N)$$
where $N$ is the cardinality of $A_P$.
\end{prop}

Let $I \subset k[(X_j), (X_j^{-1})]$ be the ideal of $V'$. 
An element $f \in I$ may be written in the form: 

$$f = \sum_{\omega \in \ze^{\oplus J}} a_\omega X^\omega $$

Where $X^\omega = \prod_{j \in J} X_j^{\omega_j}$, well defined as all $\omega_j$ but a finite number are $0$, and the set  $A_f = \{\omega \in \ze^{\oplus J} \ |\ a_\omega \neq 0\}$ is finite.

Let $f \in I$ and $s \in {\{-1,1\}}^J$, we want to separate the positive monomial of $f$ from the negative ones. Let $\omega\in A_f$. We choose an element $x \in V_s$, the number $\prod_{j \in J} f_j(x)^{\omega_j}$ has a sign depending only on $s$ and $\omega$, but not on the value of $x$. This sign is, by definition, the \nuovo{sign} of the monomial respect to $s$. We denote by $A_{f,s}^+ \subset A_f$ the set of the $\omega$'s whose monomial has positive sign respect to $s$, and by $A_{f,s}^- \subset A_f$ the set of negative ones. Now we split $f$ in its positive and negative part: 

$$f^+ = \sum_{\omega \in A_{f,s}^+} a_\omega  X^\omega \ \ \ ; \ \ \ f^- =\sum_{\omega \in A_{f,s}^-} (-a_\omega)  X^\omega $$

So $f = f^+ - f^-$.

Then we define the tropical polynomials:
$$f_0^+ = \bigoplus_{\omega \in A_f^+} 0\odot X^{\odot\omega} = \bigoplus_{\omega \in A_f^+} X^{\odot\omega}$$ 
$$f_0^- = \bigoplus_{\omega \in A_f^-} 0\odot X^{\odot\omega} = \bigoplus_{\omega \in A_f^-} X^{\odot\omega}$$

\begin{prop}
The boundary is contained in an intersection of real tropical hypersurfaces:

$$CB(V_s) \subset \bigcap_{f \in I} T_\erre(f_0^+, f_0^-)$$

In other words for all $f \in I$ a point ${(x_j)}_{j \in J} \in CB(V_s) \in \erre^J$ satisfies the relation 

$$\max_{\omega \in A_{f,s}^+} \left( \sum_{j \in J} <x,\omega>\right)  = \max_{\omega \in A_{f,s}^-} \left( \sum_{j \in J} <x,\omega> \right)$$ 

\dimo Given a point $x \in V'$ we denote the element ${(f_j(x))}_{j \in J} \in \erre^J$ by $\overline{\famil}(x)$. 
By definition of $I$ every element $x \in V_s$ verifies:

$$\forall x \in V_s : f^+(\overline{\famil}(x))  = f^-(\overline{\famil}(x))$$ 

So we have a set of equations, one for each polynomial in $I$, with coefficients in $\erre^{>0}$. 

We are interested in the set $V_{h,s}$, the image of $V_s$ under the map $D_h^J$. For every polynomial $f \in I$ we take the transformation of its positive and negative perts through $D_h$. 

$$ f_h^+ = D_h \circ f^+ \circ {(D_h^J)}^{-1} = \begin{array}{c}\ \\{\displaystyle\bigoplus}_h \\ {\scriptstyle\omega \in A_f^+}\end{array} D_h(a_\omega)\odot X^{\odot\omega}$$ 
$$f_h^- = D_h \circ f^- \circ {(D_h^J)}{-1} = \begin{array}{c}\ \\{\displaystyle\bigoplus}_h \\ {\scriptstyle\omega \in A_f^-}\end{array} D_h(-a_\omega)\odot X^{\odot\omega}$$

As $D_h$ is a semifields isomorphism we have that $V_{h,s} = \{x \in \erre^J \ |\ \forall f \in I : f_h^+(x) = f_h^-(x)\}$. 

We may put all these equation together finding equations for $W$: 
$$W = \{(x,h) \in \erre^J \times (0,+\infty) \ |\ \forall f \in I : f_h^+(x) = f_h^-(x)\}$$
Every coefficient of the polynomials $f_h^+$ and $f_h^-$ tends to $0$, the multiplicative neutral element in the semirings $S_h$, $h \geq 0$. So the functions $f_h^+$ and $f_h^-$ converges, respectively, to the functions $f_0^+$ and $f_0^-$.
This means that, as functions on $\erre^J \times (0,+\infty)$ they may be continuously extended to $\erre^J \times [0, +\infty]$ by $f_0^+$ and $f_0^-$
Hence $\overline{W}$ is contained in the set $\{(x,h) \in \erre^J \times [0,+\infty) \ |\ \forall f \in I : f_h^+(x) = f_h^-(x)\}$. So $V_{0,s}$ is contained in $\{\forall f \in I : f_0^+ = f_0^-\}$.
\end{prop}

So we have shown that the cone over the boundary is contained in an intersections of real tropical hypersurfaces, that is a polyhedral complex. In some cases it is possible to show that $CB(V_s)$ is a polyhedral subcomplex, we will do this for Teichm\"uller spaces.

We can give an estimate for the dimension of this complex, as it is contained in the cone over the boundary of the complexification. In the following section we will see that the dimension of the cone over the boundary of a complex variety is equal to the dimension of the variety.

\subsection{Compactification in the complex case}     \label{subsez:comp compl}

Let $V$ be a complex algebraic variety defined over $k\subset \ci$, and let $\famil$ be a generating family for $V$. Again we want to use Maslov dequantization on $V$.  

As before, for every $h>0$ we define the map $D_h^J$:

$$D_h^J:V' \ni x \freccia {( D_h(|f_j(x)|) )}_{j \in J} \in \erre^J$$

The image $D_h^J(V')$ will be denoted by $V_h$. Again $D_1^J$ is the $\Log$ map, so $V_1$ is the amoeba of $V$.

We define the map:

$$D:V' \times (0,+\infty) \ni (x,h) \freccia  (D_h^J(x),h)\in\erre^J \times (0,+\infty)$$

We call $W$ the image of this map. The map is proper so $W$ is closed. Let $\overline{W}$ be the closure of this set in $\erre^J \times [0,+\infty)$. As before $\overline{W}$ is simply the union of $W$ with some subset of $\erre^J \times \{0\}$, that will be called $V_0$.  

\begin{prop}
The set $V_0$ is a cone, and is equal to $CB(V')$, the cone over the boundary.

\dimo See the proof of propositions \ref{prop:cone} and \ref{prop:boundary}
\end{prop}

Let $I \subset k[(X_j), (X_j^{-1})]$ be the ideal of $V'$. An element $f \in I$ may be written in the form: 

$$f = \sum_{\omega \in \ze^{\oplus J}} a_\omega X^\omega $$

With $A_f = \{\omega \in \ze^{\oplus J} \ |\ a_\omega \neq 0\}$. We associate with $f$ a tropical polynomial $f_0$:

$$f_0 = \bigoplus_{\omega \in A_f} X^{\odot\omega}$$

We use the notation $T(I)=\bigcap_{f \in I} T(f_0)$, the intersection over all $f \in I$ of the tropical hypersurface defined by $f_0$.

\begin{prop}
Then $CB(V') \subset T(I)$.

\dimo It follows from the formula \ref{eqn:complexvar}, that holds also for abstract varieties. If $f\in I$, $f = \sum_{\omega \in A_f} a_\omega X^\omega$, then every point $x \in V_h$ satisfies:
$$\forall w \in A_f : D_h(|a_w|) \odot_h y^{\odot_h w} \leq \begin{array}{c}\ \\{\displaystyle\bigoplus}_h \\ {\scriptstyle\omega \in A_f\setminus \{w\}}\end{array} D_h(|a_\omega|) \odot_h y^{\odot_h\omega}$$
As $\lim_{h\tende 0} D_h(a_\omega) = 0$, a point $y\in V_0$ satisfies:
$$\forall w \in A_f : X^{\odot w} \leq \bigoplus_{\omega \in A_f \setminus \{w\} } X^{\odot \omega} $$ 
Hence $y\in T(f)$.
\end{prop}

We will show that $CB(V') = T(I)$ later. First we want to show how the set $T(I)$ may be described as the image of a set of valuations, as in theorem \ref{teo:bierigroves}. Let $A=k[V']=k[(f_j),(f_j^{-1})]$, the ring of coordinate of $V'$ as defined in subsection \ref{subsez:amebevarastratte}. Let $S_\erre$ be the set of real valuation of $A$ trivial over $k$. We want to show that $T(I)$ is equal to the image of the map
$$z:S_\erre \ni v \freccia {(-v(f_j))}_{j\in J} \in \erre^J$$

\begin{teo}
$z(S_\erre)=T(I)$.

\dimo When $J$ is finite the thesis follows by proposition \ref{prop:funzioni generatrici} and by \cite[thm.~2]{Be2}. The same argumentation may be extended to the general case: 

$z(S_\erre) \subset T(I)$: Let $\xi \in z(S_\erre) \subset \sferic$. There exists a valuation $v:k[V'] \freccia \erre \cup \{+\infty\}$ such that $\xi = {(-v(f_j))}_{j \in J} \in \erre^J$. This is a valuation of $k[V'] = k[(X_j), {(X_j)}^{-1}] / I$ but it may be lifted to a valuation of $k[(X_j),{(X_j)}^{-1}]$ sending the ideal $I$ in $+\infty$. Let $g \in I$, $g = \sum a_\omega X^\omega $, we have that $v(a_\omega X^\omega) = \sum_{j \in J} v(X_j)\omega_j = \sum v(f_j) \omega_j = - <\xi,\omega>$. As $v(g) = +\infty$, the value $\min_{\omega \in A_g}(v(a_\omega X^\omega))$ is assumed at least twice, so the functional $-<\xi,\omega>$ has maximum in at least two points of $A_g$, hence $\xi \in T(g)$.

$T(I) \subset z(S_\erre)$: Let $[\xi]\in D(V')$. We define a valuation in the following way:  
$$v: k[(X_j),{(X_j)}^{-1}] \ni f \freccia v(f) = \inf_{\omega \in A_f} - <\xi,\omega> \in \erre\cup\{+\infty\}$$
We denote by $S\subset k[(X_j),{(X_j)}^{-1}]$ the multiplicative subgroup generated by the symbols $X_j$. The ideal of $V'$, $I$, does not intersect $S$, moreover if $f \in I$, we know that the minimum value of $-<\xi,\omega>$ over $A_f$ is achieved at least twice, hence there are at least two monomials of $f$ with minimal valuation, say $X^{\omega'}, X^{\omega''}$. This guarantees that the valuation of $f-X^{\omega'}$ is equal to $v(f)$. If $s \in S$, and if $v(f) = v(s)$, we have that $v(f-s) \geq \min( v(X^{\omega'}-s), v(f-X^{\omega'}) \geq v(f)$. Hence the hypothesis of \cite[coroll.~1]{Be1} holds. This theorem, together with the note that follows it, guarantees the existence of a valuation $v'$ of $k[(X_j),{(X_j)}^{-1}]$ such that: 
$\forall f\in S: v'(f) = v(f)$ and $\forall f\in I: v'(f)=+\infty$. By lemma \ref{lemma:quot val}, this valuation defines a valuation $v''$ of $k[V']$ such that $z(v'')=\xi$. 
\end{teo}

Now we want to show that $CB(V') = T(I)$. For an hypersurface this is a particular case of \cite[coroll. 6.4]{Mi}.

\begin{prop}
If $V \subset \ci^n$ is the hypersurface defined by the polynomial $f% = \sum_{\omega \in \ze^n} a_\omega X^\omega 
$, and the generating family $\famil$ is the family of coordinates $\{X_1 \dots X_n\}$, then $CB(V') = %$ is equal to the tropical hypersurface defined by the polynomial $f_0 = \bigoplus_{\omega \in A_f} X^{\odot\omega}$.
T(f_0)$.

\dimo Let $\cappa$ be the field of Puiseux series in $t$ with complex coefficients. The coefficients of $f$ may be viewed as constant series, becoming elements of $\cappa$ whose valuation is $0$. If we apply \cite[coroll. 6.4]{Mi} to the polynomial $f$, it states that the sets $V_h$ converge to $V_0$ in the Hausdorff metric when $h$ tends to $0$. This fact implies the conclusion.
\end{prop}

When $\famil$ is finite it follows from a result from \cite{Be2} and one from \cite{BG}.

\begin{prop}
If $\famil$ is finite, then $CB(V')=T(I)$.

\dimo We only need to show that $T(I) \subset CB(V')$. By proposition \ref{prop:funzioni generatrici} we may suppose that $V \subset \ci^n$ and that $\famil=\{X_1 \dots X_n\}$. In \cite{Be2}, in step (4) of the proof of thm. 2, it is shown that all points of $T(I)$ with rational coordinates are in $CB(V')$. Then in \cite{BG} it is shown that rational points are dense in $T(I)$.
\end{prop}

The same result is true when $\famil$ is infinite, and will be proved in section \ref{sez:succ e  val}.

If $J$ is finite the dimension of $B(V')$ may be computed:

\begin{teo}  \label{teo:dimensione}
If $J$ is finite and if $V$ is irreducible and $\dim_\ci(V) = m$, then $B(V')$ is a finite union of $(m-1)$-dimensional subspheres.

\dimo See \cite[thm.~3]{Be2}.
\end{teo}

\section{Algebraic structure on Teichm\"uller spaces}  

To apply the construction presented above to Teichm\"uller spaces we need an ``algebraic structure'' on them. This topic was treated in \cite{MS1}, using some theory developed in \cite{CS1}. In this section we sketch the construction of the ``algebraic structure'' and we fix some notations that will be useful in the following. 

Note that the variety $X^{par}_\erre(\Gamma)$ that we will define is a real variety defined over $\qu$, and it is well defined up to polynomial isomorphisms as it is the real part of a complex variety.  

\subsection{Variety of characters}  \label{subsez:caratteri}

Let $\Gamma$ be a finitely generated group, and let $R(\Gamma)$ be the \nuovo{variety of representations} of $\Gamma$ in $\sld$, $R(\Gamma) = \{\rho:\Gamma \freccia \sld\}$, endowed with its structure of complex affine variety defined over $\qu$. We denote by $\tau_\gamma \in \qu[R(\Gamma)]$ the functions $\tau_\gamma:R(\Gamma) \ni \rho \freccia \tr(\rho(\gamma)) \in \ci$.

We define the \nuovo{character} of a representation $\rho \in R(\Gamma)$ to be the function $\chi_\rho:\Gamma \ni \gamma \freccia \tr(\rho(\gamma)) \in \ci$. We will denote by $X(\Gamma)$ the set of all character of representations in $R(\Gamma)$. This set too may be endowed with a structure of affine complex variety defined over $\qu$ (see \cite{CS1}), a structure that turns the natural map $t:R(\Gamma) \ni \rho \freccia \chi_\rho \in X(\Gamma)$ in a regular map over $\qu$. This variety will be called \nuovo{variety of characters}.   

For $\gamma\in\Gamma$, the value of $\tau_\gamma$ only depends on the character of its argument, hence it induces a function $I_\gamma:X(\Gamma) \ni \chi_\rho \freccia \tau_\gamma(\rho) = \tr(\rho(\gamma))$. The functions $I_\gamma$ will be called \nuovo{trace functions}. The trace functions belong to the coordinate ring $\qu[X(\Gamma)]$. Moreover there exists a finite set $A\subset \Gamma$ such that the family ${\{I_\gamma\}}_{\gamma \in A}$ generates the ring $\qu[X(\Gamma)]$ as a $\qu$-algebra (see \cite{CS1}).

The function $I_\gamma$ only depends on the conjugation class $c$ of $\gamma$, hence we can write it as $I_c$. We denote by $\class$ the set of conjugation classes of elements in $\Gamma$, and we denote the family of all trace functions by $\gfamil = \{I_c\}_{c\in\class}$. This is a countable family, and it is an example of what will be called a generating family. 

We want to focus on representations of $\Gamma$ in $SL_2(\erre)$, corresponding to points of the real part of $R(\Gamma)$, $R_\erre(\Gamma)$. The characters of these representations are in the real part of $X(\Gamma)$, $X_\erre(\Gamma)$ because the map $t$ is a regular map over the field $\qu$, but the image $t(R_\erre(\Gamma)) \subset X_\erre(\Gamma)$ is generally a strict subset of $X_\erre(\Gamma)$. %For example the character of a representation whose image is contained in $SU_2(\ci)$ is in $X_\erre(\Gamma)$.

\subsection{Fundamental groups of surfaces}

Let $S$ be a compact connected orientable surface of genus $g$ with $k$ boundary components and let $\Gamma=\pi_1(S)$. If $k>0$, for every boundary component we may choose an element $\beta_i\in [\Gamma, \Gamma]$ that has the same free homotopy class. 

Suppose, in the following, that $\chi(S) < 0$. So we can put on $\interno{S}$ a complete hyperbolic structure of finite volume. The holonomy representation of this structure is a discrete and faithful representation $\rho:\Gamma \freccia PSL_2(\erre)$, such that for every component of $\partial S$ the corresponding element $\beta_i$ is parabolic, i.e. $|\tr(\rho(\beta_i))|=2$. The Teichm\"uller space of $S$ may be defined as the set of conjugation classes of all representations satisfying these properties. This space depends only on $g$ and $k$, and will be denoted by $\tau_g^k$.

The sets $R^{par}(\Gamma) = \{ \rho \in R(\Gamma) \ |\ \forall i: \tr(\rho(\beta_i)) = \pm 2 \}$ and $X^{par}(\Gamma) = \{ \chi \in X(\Gamma) \ |\ \forall i: I_{\beta_i}(\chi) = \pm 2 \}$ are affine subvarieties of $R(\Gamma)$ and $X(\Gamma)$ respectively, both defined over $\qu$. Again we focus on real parts, $R^{par}_\erre(\Gamma) = R^{par}(\Gamma) \cap R_\erre(\Gamma)$ and $X^{par}_\erre(\Gamma) = X^{par}(\Gamma) \cap X_\erre(\Gamma)$.

We want to identify the Teichm\"uller space of $S$ with a connected component of $X^{par}_\erre(\Gamma)$. The keypoint is a result from Thurston's stating that all representation whose conjugation class belongs to $\tau_g^k$ may be lifted to representations taking values in $SL_2(\erre)$.

We define $DR(\Gamma)\subset R(\Gamma)$ as the subset of discrete faithful representations, $DX(\Gamma) = t(DR(\Gamma)) \subset X(\Gamma)$ and their ``real parts'', $DR_\erre(\Gamma) = DR(\Gamma) \cap R_\erre(\Gamma)$, and $DX_\erre(\Gamma) = DX(\Gamma) \cap X_\erre(\Gamma)$.

The conjugation class of a representation in $DR(\Gamma)$ may be identified with its character, as representations in $DR(\Gamma)$ are irreducible (see \cite[lemma~3.1.3]{MS1}) and by \cite[prop.~1.5.2]{CS1} if an irreducible representation has the same character as an other representation, then they are conjugated.

The subset of $R(\Gamma)$ that is closer to $\tau_g^k$ is $DR^{par}_\erre(\Gamma) = DR(\Gamma) \cap R^{par}_\erre(\Gamma)$. By \cite[prop.~3.1.4]{MS1} the set of characters of these representation is equal to the set $DX^{par}_\erre(\Gamma) = DX(\Gamma) \cap X^{par}_\erre(\Gamma)$. 

The subset $DR^{par}_\erre(\Gamma)$ may be seen as the set of representations in $SL_2(\erre)$ that composed with the quotient in $PSL_2(\erre)$ are the holonomy of a hyperbolic structure on $S$ (they remain faithful as $\Gamma$ is torsion-free). The set $DX^{par}_\erre(\Gamma)$ may be seen as the set of conjugation classes of these representations. This set is a union of connected components (with reference to the classical topology) of $X^{par}_\erre(\Gamma)$ (see \cite[prop.~3.1.8]{MS1}). We want to identify the Teichm\"uller space of $S$ with one of these components.

We consider the following action of $H = \homom(\Gamma,\ze_2)$ on $DX^{par}_\erre(\Gamma)$: 
$H \times DX^{par}_\erre(\Gamma) \ni (h,\chi_\rho) \freccia \chi_{\rho'} \in DX^{par}_\erre(\Gamma)$, where $\rho'(\gamma) = {(-1)}^{h(\gamma)} \rho(\gamma)$ and $\chi_{\rho'}(\gamma) = {(-1)}^{h(\gamma)} \chi_\rho(\gamma)$. Two elements have the same orbit if and only if they induce the same conjugation class of representation in the quotient $PSL_2(\erre)$, hence the quotient by this action may be identified with $\tau_g^k$.  

\begin{lemma}  \label{lemma:segno costante}
Trace functions $I_\gamma$ never vanish on $DX^{par}_\erre(\Gamma)$.

\dimo Suppose $I_\gamma(\chi_\rho)=0$. This means $\rho(\gamma) \in \sld$ and $\tr(\rho(\gamma)) = 0$, hence the Jordan form of $\rho(\gamma)$ is 
$\begin{pmatrice}
i &  0 \\
0 & -i
\end{pmatrice}$. The matrix $\rho(\Gamma)$ has finite order, but $\Gamma$ is a torsion-free group, hence $\rho$ may not be faithful. 
\end{lemma} 

Let $h \in H$, $h \neq \ident$. There exists $\gamma \in \Gamma$ such that $h(\gamma) = -1$, hence $\forall \chi: I_\gamma(\chi) = - I_\gamma(h\chi)$. By lemma \ref{lemma:segno costante}, if $\chi$ is a character,$\chi$ and $h\chi$ belong to different connected components of $DX^{par}_\erre(\Gamma)$. In other words $H$ acts freely on the set of connected components of $DX^{par}_\erre(\Gamma)$.   

The identification with $DX^{par}_\erre(\Gamma) / H$ induces a topology on $\tau_g^k$, such that it is homeomorphic to an euclidean ball. So the quotient is connected, and this implies that the action of $H$ on the set of connected components of $DX^{par}_\erre(\Gamma)$ is transitive. Hence the quotient $DX^{par}_\erre(\Gamma) / H$ may be identified with any of the connected components of $DX^{par}_\erre(\Gamma)$, each of which is a connected component of $X^{par}_\erre(\Gamma)$.

\subsection{Length and trace functions}        \label{subsez:length and trace functions}

We have seen that the Teichm\"uller space of a surface $S$ may be identified with a connected component of $X^{par}_\erre(\Gamma)$ (where $\Gamma = \pi_1(S)$), the real part of $X^{par}(\Gamma)$. We have also defined the family of functions $\gfamil = \{I_c\}_{c\in \class}$. If two points of $DX^{par}_\erre(\Gamma)$ has the same orbit for the action of $H$, the values of functions $I_c$ on these points coincide up to sign. So $|I_c|$ is a well defined function on $\tau_g^k$, and will be called, again, a trace function. Trace functions are closely related to length functions: $\ell_c([h]) = \inf_{\alpha} l_h(\alpha)$, where $h$ is an hyperbolic metric on $S$, $[h]$ is the corresponding elements of $\tau_g^k$, $l_h$ is the function that send a curve in its $h$-length, and the $\inf$ is taken on the set of all closed curves whose free homotopy class is $c$. The relation between trace functions and length function is given by $|I_c([h])| = 2 \cosh(\frac{1}{2}\ell_c([h]))$. This implies that $|I_c([h])|\geq 2$.

\section{Compactification of Teichm\"uller spaces}

\subsection{Thurston boundary}

Let $S$ be a compact connected surface with genus $g$ and with $k$ boundary components, and let $\Gamma = \pi_1(S)$. The set $DX^{par}_\erre(\Gamma)$ is a union of connected components of $X^{par}_\erre(\Gamma)$. We choose $\gfamil = \{I_c\}_{c\in \class}$, as a generating family for $X^{par}_\erre(\Gamma)$, the set $J = \class$ being countable. As no elements of $\gfamil$ vanish on $DX^{par}_\erre(\Gamma)$, is well defined the compactification, as in subsection \ref{subsez:varreali}, with $\overline{DX^{par}_\erre(\Gamma)} = DX^{par}_\erre(\Gamma) \cup B(DX^{par}_\erre(\Gamma))$.

As the action of $H = \homom(\Gamma,\ze_2)$ commutes with the $\Log$ map, we may extend the identification between two connected components induced by an element of $H$ to an identification between the compactifications. This way we may construct a compactification of the Teichm\"uller space: $\overline{\tau_g^k} = \tau_g^k \cup B(\tau_g^k)$. The boundary $B(\tau_g^k)$ may be seen as a subset of $\sferic$. This is because the boundary of a connected component of $DX^{par}_\erre(\Gamma)$, as a subset of $\sferic$, does not depend on the chosen component.

We may define the $\Log$ map on $\tau_g^k$ as:
$$\Log:\tau_g^k \ni h \freccia {(\log(|I_c(h)|))}_{c \in \class} \in \erre^{\class}$$
Hence we have: $B(\tau_g^k) = \{ b \in \sferic \ |\ \exists \{x_n\} \subset \tau_g^k: x_n \tende \infty$ in $\overline{\tau_g^k}$ and $\pi(\Log(x_n)) \tende b$ in $\sferic \}$.

\begin{prop}      \label{prop:Thurston}
This compactification of Teichm\"uller spaces is the same as the one described by Thurston, with a correspondence that sends the projective class of a measured foliation $f$ in the point ${[I(f,c)]}_{c \in \class} \in \sferic$, $I(f,c)$ being the measure of $c$ with reference to $f$.

\dimo If $(x_n) \subset \tau_g^k$ is a sequence converging to $f$ in Thurston compactification, then ${[I(f,c)]}_{c \in \class}$ is the limit of the sequence ${[\ell_c(x_n)]}_{c\in \class}$, but this is equal to the limit of the sequence ${[I_c(x_n)]}_{c \in \class}$ (see subsection \ref{subsez:length and trace functions}). 
\end{prop}

Let $A \subset \class$, and suppose the family $\famil_A = \{I_c \ |\ c \in A\}$ to be a generating family, or a finite family. This family defines a compactification of $DX^{par}_\erre(\Gamma)$, and, as before, a compactification of $\tau_g^k$. Suppose, from now on, the map $\overline{\famil}$ to be proper on $DX^{par}_\erre(\Gamma)$. We will denote by $B_A(\tau_g^k)$ the boundary constructed using the family $\famil_A$, and by $CB_A(\tau_g^k)$ the cone over this set.  

Let $A \subset B \subset \class$, we denote the canonical projection by $\pi_{B,A}:\erre^B \sfreccia \erre^A$. 

\begin{lemma}
$\pi_{B,A}(CB_B(\tau_g^k)) \subset CB_A(\tau_g^k)$.

\dimo We will use the description of $CB_B(\tau_g^k)$ and $CB_A(\tau_g^k)$ in terms of valuations. Let $x \in CB_B(\tau_g^k)$. If $\pi_{B,A}(x) = 0$ there is nothing to prove. Else there exists a valuation $v$ of the ring generated by $\famil_B$ such that $z(v) = x$. We denote by $w$ the valuation $v$ restricted to the ring generated by $\famil_A$. $w \in S_\erre'$ as $\pi_{B,A}(x) \neq 0$. Hence $\pi_{B,A}(x) = z(w) \in CB_A(\tau_g^k)$.
\end{lemma}

Let $A \subset B \subset \class$ be as before. We say $A$ and $B$ to be \nuovo{compatible} if $\pi_{B,A}^{-1}(0) \cap CB_B(\tau_g^k) = \{0\}$. With this hypothesis the restricted projection $\pi_{B,A|CB_B(\tau_g^k)}$ induces a map between the spherical quotients: $\overline{\pi}_{B,A}:B_B(\tau_g^k) \freccia B_A(\tau_g^k)$. 

The compatibility relation is transitive. If $A \subset B \subset C$ we have that $\overline{\pi}_{C,A} = \overline{\pi}_{C,B} \circ \overline{\pi}_{B,A}$. 

If $A \subset B \subset \class$ are compatible family we may define a map $p:\tau_g^k \cup B_B(\tau_g^k) \freccia \tau_g^k \cup B_A(\tau_g^k)$, defined by the identity on $\tau_g^k$ and by the map $\overline{\pi}_{B,A}$ over $B_B(\tau_g^k)$.

\begin{prop}
The map $p$ defined above is continuous.

\dimo Continuity over $\tau_g^k$ is obvious. Now we take a point $x\in B_B(\tau_g^k)$ and a sequence  $\{x_n\} \subset \tau_g^k \cup B_B(\tau_g^k)$ converging to $x$. We only need to show that every subsequence $\{y_n\}$ has a subsequence $\{y_{n_k}\}$ such that $p(y_{n_k})$ tends to $p(x)$. If $\{y_n\}$ is infinitely often in $B_B(\tau_g^k)$, we may take a subsequence contained in the boundary, and we may conclude by continuity of $\overline{\pi}_{B,A}$. Else $\{y_n\}$ is infinitely often in $\tau_g^k$, so we may extract a subsequence $\{y_{n_k}\}$ contained in $\tau_g^k$ and quasi-valuating. Let $v$ be a valuation compatible with $\{y_{n_k}\}$. We know that $\{y_{n_k}\}$ converges to ${[-v(I_c)]}_{c \in B}$ in $\tau_g^k \cup B_B(\tau_g^k)$, and that it converges to ${[-v(I_c)]}_{c \in A}$ in $\tau_g^k \cup B_A(\tau_g^k)$. The former point is $x$, and the latter is $p(x)$.  
\end{prop}

\begin{corol}
The map $p$ is surjective (so also the map $\overline{\pi}_{B,A}$ is).

\dimo $p$ is a continuous map, from a compactum to an Hausdorff space, so it is a closed map. The image is dense (as it contains $\tau_g^k$), so the map is surjective.
\end{corol}

Let $A\subset\class$. We define $A$ to be \nuovo{equivalent} to $\class$ if $A$ is compatible with $\class$ and $\overline{\pi}_{\class,A}$ is injective. With this hypothesis the map $p$ between the compactification defined by $\class$ and the one defined by $A$ is an homeomorphism. This is the proof of the following proposition:

\begin{prop}
The Thurston boundary for Teichm\"uller space $\tau_g^k$ may be constructed starting by any family $\famil_A$ with $A$ equivalent to $\class$.
\end{prop}

\subsection{PL-structure on the boundary}          \label{subsez:PL struc}

Now it is useful to find some finite subset of $\class$ that are equivalent to $\class$. We start finding some compatible subsets.

\begin{prop}
Let $A\subset \class$. If $\famil_A$ is a generating family for $X^{par}(\Gamma)$, then $A$ is compatible with $\class$.

\dimo By contradiction, we suppose the existence of a non zero element $x \in CB(\tau_g^k)$ such that $\pi_{\class,A}(x) = 0$. There exists a valuation $v: k[V'] \freccia \erre\cup \{+\infty\}$ such that $\forall c \in A : v(I_c) = 0$, and $\exists \overline{c} \in \class : v(I_{\overline{c}}) < 0$. This is impossible as $I_{\overline{c}}$ is a polynomial in the elements of $\famil_A$, so $v(I_{\overline{C}}) \geq \max_{c\in A}\{v(I_c)\} \geq 0$.
\end{prop}

\begin{prop}
Suppose that $A \subset \class$ is a system of free homotopy classes of curves that \nuovo{fill up}, i.e. every free homotopy class of closed curves on the surface has non zero intersection form with at least one of those curves. Then $A$ is compatible with $\class$

\dimo Let $x\in B(\tau_g^k)$. Let $f$ and $v$ be respectively a measured foliations and a valuation associated to $x$. From proposition \ref{prop:Thurston} we know that if $I(f,c) \neq 0$, then $v(I_c) \neq 0$. 

If $A$ is a system that fills up we have that for all $x\in B(\tau_g^k)$ there exists a $c\in A$ such that $I(f,c) \neq 0$, and this implies $v(I_c)\neq 0$. So $A$ is compatible.

Moreover we have shown that is $A$ fills up the point ${[v_x(I_\gamma)]}_{\gamma \in A} \in \sferic$ is equal to the point ${[I(f_x,\gamma)]}_{\gamma \in A}\in \sferic$. This will be useful in the proof of the following proposition.   
\end{prop}

\begin{prop}
There exists a finite set $A \subset \class$, consisting on $9g-9+3b$ free homotopy classes of simple curves, that is equivalent to $\class$.

\dimo There exists $3g-3+b$ simple curves on $S$ (denoted by $K_1 \dots K_{3g-3+b}$) that decompose it in $2g-2+b$ pair-of-pants, $b$ of them containing a boundary component of $S$. The curve $K_i$ is the common boundary of two pair-of-pants whose union will be denoted by $P_i$. We denote by $K_i', K_i''$ the two simple closed curves in $P_i$ defined by Thurston in the classification of measured foliation (See \cite{FLP}). Let $A$ be the set of the free homotopy classes of the curves $K_i$, $K_i'$ and $K_i''$. This set fills up, so $A$ is compatible with class. To show that it is equivalent we need to show that the map $\overline{\pi}_{\class,A}$ is injective. This follows from Thurston's classification of measured foliation (see also the proof of the previous proposition). 
\end{prop}

Now we fix a connected component $C \subset DX^{par}_\erre(\Gamma)$, and we denote by $X$ the intersection of all subvarieties of $X^{par}(\Gamma)$ containing $C$. The variety $X$ is a complex irreducible variety satisfying ${(X_\erre)}_\ci = X$. Let $m$ be the (complex) dimension of $X$, equal to the real dimension of $C$. Let $A \subset \class$ a finite set equivalent to $\class$. Let $X'^{comp}$ be the compactification of $X'$ with reference to the family ${\{I_{c|X}\}}_{c \in A}$. The closure of $C$ in ${X'}^{comp}$ may be identified with the compactification of Teichm\"uller space ${\tau_g^k}^{comp} = \tau_g^k \cup B_A(\tau_g^k)$. The boundary $B_A(\tau_g^k)$ is homeomorphic to a sphere $S^{m-1}$, and is contained in $B_A(X')$. This latter set is, by theorem \ref{teo:dimensione}, a finite union of $(m-1)$-dimensional subspheres, hence it may be seen as a $(m-1)$-dimensional polyhedral complex.

\begin{lemma}
Let $M$ be a compact topological $n$-manifold without boundary, let $P$ be a $n$-dimensional polyhedral complex, and let $i:M \ifreccia P$ be a injective continuous map. Then for every $n$-dimensional face $F \subset P$ the intersection $i(M) \cap \interno{F}$ is empty or the whole $\interno{F}$. Hence the image $i(M)$ is the union of all faces of $P$ whose interior intersects $i(M)$.

\dimo The set $D = \{x \in M \ |\ i(x) $ is in the interior of an $n$-dimensional face of $P \}$ is a dense open subset of $M$, as it is impossible for a neighborhood of a point to be sent injectively in the $(n-1)$-skeleton. We fix an $n$-dimensional face $F \subset P$, and let $F' = i^{-1}(F) \cap D$. We want to say that $i(F')$ is the whole $\interno{F}$ or it is empty. It is enough to say that it is open and closed. It is open as the map $i_{|F'}:F' \freccia \interno{F}$ is an injective map from a $n$-manifold to a space homeomorphic to $\erre^n$, so, by the invariance of domain, it is an open map. To show that it is closed we take a point $x$ in the closure. There exists some $y \in M : i(y)=x$ as $M$ is compact, so the image is closed. The point $y$ is necessarily in $D$ (by definition of $D$), so $y\in F'$. Hence $x \in i(F')$.    
\end{lemma}

\begin{prop}
The subset $B_A(\tau_g^k)$ is a polyhedral subcomplex of $B_A(X')$.

\dimo It is enough to apply previous lemma to the canonical immersion $i:B_A(\tau_g^k) \ifreccia B_A(X')$. 
\end{prop}

So if we fix a finite set $A\subset \class$ equivalent to $\class$ we have that $B_A(\tau_g^k)$ is a polyhedral complex.

\begin{prop}
Let $A,B \subset \class$ finite subsets equivalent to $\class$. Then $B_A(\tau_g^k)$ and $B_B(\tau_g^k)$ are PL-homeomorphic polyhedral complexes.

\dimo We have that $A\cup B$ is finite and equivalent to $\class$. The two maps $\pi_{A\cup B, A}$ and $\pi_{A\cup B, B}$ are restriction of linear maps, so they are PL-maps. Hence the maps they induce on the spherical quotient:  $\overline{\pi}_{A\cup B, A}$ and $\overline{\pi}_{A\cup B, B}$ are PL-homeomorphisms, so ${(\overline{\pi}_{A\cup B, A})}^{-1} \circ \overline{\pi}_{A\cup B, B}: B_A(\tau_g^k) \freccia B_B(\tau_g^k)$ is a PL-homeomorphism.
\end{prop}

This way we have constructed a PL-structure on the boundary independently from the chosen finite set $A \subset \class$. Now we want to show that this PL-structure is PL-homeomorphic to the one defined by Thurston using train tracks. See \cite{Pap} for definitions and details. 

An admissible train track on $S$, is a graph embedded in $S$ satisfying certain conditions. A measure on a train track is a function from the set of the edge in $\erre^{> 0}$, satisfying, again, certain conditions. There exists an enlargement operation associating to every measured admissible train track a measured foliation $f$ on $S$, and so a point of $CB(\tau_g^k)$. 

If $\tau$ is a fixed train track with $n$ edges, every measure on $\tau$ may be identified with a point of ${\erre^{> 0}}^n$, and the subset of all these point is a polyhedral conic subset, that will be denoted by $C_\tau$. The enlargement operation defines a map $\phi_\tau: C_\tau \freccia CB(\tau_g^k)$. If every connected component of $S \setminus \tau$ is a triangle, the image $\phi_\tau(C_\tau)=V_\tau$ is an open subset of $CB(\tau_g^k)$. Moreover the map $\phi_\tau$ is an homeomorphism with its image. 

The union of the open sets $V_\tau$ is the whole $CB(\tau_g^k)$. So every $V_\tau$ is identified with $C_\tau$ in such a way that the changes of charts are piecewise linear. This way we have described a PL-atlas for $CB(\tau_g^k)$, the PL-structure defined by Thurston.

We want to show that the identity map is a PL map if we endow the domain with the Thurston's PL-structure, and the codomain with the structure defined above. This implies that the two structure are equal.

To show this we need to prove that the maps $\phi_\tau$ are PL if we endow the codomain with the structure defined above. We choose a finite family $A \subset \class$ equivalent to class, so that the cone $CB_A(\tau_g^k)$ is a subset of $\erre^A$. %The set $C_\tau$ is also a subset of $\erre^n$, $n$ being the number of edges of $\tau$. 
The map $\phi_\tau: C_\tau \freccia CB_A(\tau_g^k) \subset \erre^A$ is made up of coordinates, for each element $c\in A$ the corresponding coordinate of $\phi_\tau$ is the function that associate to a measure $\mu$ on $\tau$ the number ${(I(f,c))}$, where $f$ is the foliation constructed by the enlargement of the measured train track $(\tau,\mu)$.
 
For all $c \in A$ it is easy to see that the corresponding coordinate is PL. We choose a curve $\gamma \in c$ such that $\gamma$ does not contain any vertex of $\tau$ and such that it intersects any edge transversely. Now we define the function $p_\gamma: C_\tau \freccia \erre$ as the sum of the measures of all the edges intersected by $\gamma$, counted with multiplicity. This function is the restriction of a linear function with positive integer coefficients. The coordinate of $\phi_\tau$ corresponding to $c$ is simply the minimum of all the function $p_\gamma$, and locally the minimum may be taken over a finite numbers of linear functions, so the result is a PL function.

\subsection{Compactification using an hypersurface}

It may be useful to find a set $A \subset \class$ equivalent to $\class$ of minimal cardinality. For $\tau_g^k$ it may be found such a set with $6g-5+2b$ elements, just one more than the dimension. What we need is a set of curves that fill up and such that the map ${(I(\cdot,c))}_{c \in A}$ from the set of Whithead class of measured foliation in $\erre^A$ is injective. If $k=0$ such a set is described in \cite{Ha2}, else it is described in \cite{Ha1}. In the following $A$ will denote a set with these properties.

Let $X$ be, as in subsection \ref{subsez:PL struc}, the intersection of all subvarieties of $X^{par}(\Gamma)$ containing a fixed connected component $C$ of $DX^{par}_\erre(\Gamma)$. The variety $X$ is a real variety of dimension $6g-6+2b$. Let $\overline{\famil_A}$ be the map:

$$\overline{\famil_A}:X \ni x \freccia {(I_c(x))}_{c \in A} \in \erre^A \simeq \erre^{6g-5+2b}$$

We denote by $W$ the closure of the image of $\overline{\famil_A}$. This is an irreducible variety, and its dimension is lesser or equal than the dimension of $X$ $(6g-6+2b)$. The boundary of $C$ with reference to the family $A$ is, by definition (see subsection \ref{subsez:non generating}) contained in the boundary of $W$ with reference to the family of coordinates $(X_1 \dots X_{6g-5+2b})$. So the boundary of $C$ has dimension greater or equal to $6g-7+2b$, and this implies that the dimension of $W$ is exactly $6g-6+2b$. So $W$ is an hypersurface in $\erre^{6g-5+2b}$, and the cone over its boundary is contained in a tropical hypersurface. 

This way we identify the cone over the boundary of $\tau_g^k$ with a subpolyhedron of a tropical hypersurface in $\erre^{6g-5+2b}$. 

The map $\overline{\famil_A}$ is proper on $C$. If $k>0$ and $A$ is the set described in \cite{Ha1} this map is also injective. This may be shown using the fact that the map 
$$\tau_g^k \ni x \freccia {(\ell_x(c))}_{c\in A} \in \erre^A$$
is injective (see \cite{Ha1}), and that $\famil_A$ is the composition of this map with $\cosh$.  

So we have a homeomorphism from $\tau_g^k$ to a closed subset of $W$.

Unluckily we currently don't know the equation for $W$, nor for the associated real tropical hypersurface. This piece of information could be very useful. For example, if we knew that the equation has few monomials, we could say that the tropical hypersurface. As this surface is a cone, this would imply that it is a disk, hence that it is a cone over a sphere. This would have two nice consequences: we could say that the cone over the boundary is the whole tropical hypersurface, not only a subpolyhedron, and we could have an alternative proof that the boundary for the Teichm\"uller space is a sphere.

\section{Appendix: Valuations and sequences}    \label{sez:succ e  val}

The main target of this section is to prove that $CB(V') = T(I)$, as stated in subsection \ref{subsez:comp compl}, in the case of a infinite family $\famil$. We make this using a modification of Morgan and Shalen theory of valuating sequences, that is also a nice way for interpreting the description of boundary points as valuations.

\subsection{Ordered groups}  \label{subsez:gruppi ord}

In this subsection we fix some notations about ordered groups. By an \nuovo{ordered group} we mean an Abelian group $(\Lambda,+)$ equipped with a total order relation $<$ satisfying the property: $\forall a,b : a < b$ and $c \leq d$ $\Rightarrow$ $a+c < b+d$. We define $|\lambda| = \max\{\lambda,-\lambda\}$. By a \nuovo{closed interval} with extremes $\lambda_1,\lambda_2 \in \Lambda$ we mean the set $[\lambda_1,\lambda_2] = \{\lambda \in \Lambda \ |\ \lambda_1 \leq \lambda \leq \lambda_2 \}$. A convex subset is a set $I \subset \Lambda$ such that $\forall \lambda_1, \lambda_2 \in I: [\lambda_1,\lambda_2] \subset I$.

Every subgroup of an ordered group is viewed as an ordered group with the restricted order relation. If $\Lambda' \subset \Lambda$ is a convex subgroup, the order relation of $\Lambda$ induces an order relation on the quotient $\Lambda / \Lambda'$.

The set of convex subgroups of $\Lambda$ is totally ordered by inclusion. If this set is finite the group is said to be of \nuovo{finite rank}, and the \nuovo{rank} of the group is defined to be the cardinality of the set of its non trivial convex subgroups. 

The basic example of a ordered group of rank $n$ is $\erre^n$ with the lexicographic order. Every subgroup of rank $n$ is isomorphic to some subgroup of $\erre^n$. If the rank is $1$ the immersion of the group in $\erre$ is uniquely determined up to scalar multiplication. 

Let $\Lambda$ be a ordered group of finite rank $r$. We denote its convex subgroups by $\{0\} = \Lambda_0 \subsetneq \dots \subsetneq \Lambda_r = \Lambda$. The \nuovo{height} of an element $\lambda \in \Lambda$ is defined as the $j\in \{1 \dots r\}$ such that $\lambda \in \Lambda_j \setminus \Lambda_{j-1}$. The group $\Lambda_j / \Lambda_{j-1}$ has rank $1$, hence there exists an immersion $i_j: \Lambda_j / \Lambda_{j-1} \ifreccia \erre$, unique up to a scalar factor.  

We endow the set $[-\infty,+\infty] = \erre \cup \{\pm \infty\}$ with extended operations such that $\frac{1}{\pm\infty} = 0$, $a\cdot (\pm\infty) = (\pm\segno(a))\infty$. 

We introduce a division operation $\Lambda \times (\Lambda\setminus\{0\}) \ni (\lambda,\lambda') \freccia \frac{\lambda}{\lambda'} \in \erre\cup\{\pm\infty\}$, in the following way: if $\lambda$ and $\lambda'$ have the same height $j$ we define $\frac{\lambda}{\lambda'} = \frac{i_j(\lambda)}{i_j(\lambda')}$, if $\lambda$ has an height less than the height of $\lambda'$ we define $\frac{\lambda}{\lambda'} = 0$, if it is greater we define $\frac{\lambda}{\lambda'} = \segno(\lambda)\segno(\lambda')\infty$. 

With this definition the following properties hold: $\forall \lambda,\lambda' : (\frac{\lambda}{\lambda'}) \neq 0 \Rightarrow {(\frac{\lambda}{\lambda'})}^{-1} = \frac{\lambda'}{\lambda}$ and  $\forall r,s\in \ze$ with $\segno(s)= -\segno(\lambda')$, we have $\frac{\lambda}{\lambda'} < \frac{r}{s}$ $\Leftrightarrow$ $ r\lambda' < s\lambda$.

\subsection{Ring valuations}

In this section we give a concise introduction to ring valuation, and then we show how one can use them to make a correspondence from a subset of the set of finite rank valuation of a field to the set of rank one valuation of a subring.

\begin{defin}
Let $k$ be a field, $A \supset k$ a $k$-algebra and $(\Lambda,+,<)$ an ordered group. By a \nuovo{valuation} of $A$ over $k$ we mean a map $v:A \freccia \Lambda \cup \{+\infty\}$ such that: 
\begin{enumerate}
	\item $v(xy) = v(x) + v(y)$.
	\item $v(x + y) \geq \min(v(x),v(y))$.
	\item $v(0) = \infty$, $v(k^*) = 0$.
\end{enumerate}

A valuation is said to be \nuovo{trivial} if $v(A) \subset \{0,+\infty\}$.
\end{defin}

The set $v^{-1}(+\infty)$ is a prime ideal of $A$. Hence when $A=\effe$ is a field we have $v(\effe^*) \subset \Lambda$.

The image $v(A)\setminus\{+\infty\}$ is a subgroup of $\Lambda$, and will be called $\Lambda^v$. If $\Lambda^v$ is of finite rank $v$ is said to be of \nuovo{finite rank} and the \nuovo{rank} of $v$ is defined as the rank of $\Lambda^v$. 
 
\begin{defin}
Two valuations $v:A \freccia \Lambda \cup \{+\infty\}$ and $w:A \freccia \Lambda'\cup\{+\infty\}$ are said to be \nuovo{equivalent} if $v^{-1}(+\infty) = w^{-1}(+\infty)$ and there exists a ordered group isomorphism $\varphi:\Lambda^v \bfreccia \Lambda^w$ such that $\forall x \in A\setminus v^{-1}(+\infty) : \varphi(v(x)) = v'(x)$.
\end{defin}

When we are interested to valuation up to equivalence, we may suppose that they are surjective ($\Lambda^v = \Lambda$). Else if we want to compare valuations whose values belongs to the same group (as, for example, real valuations, $\Lambda = \erre$) we may not suppose surjectivity.

\begin{defin}
If $v$ is a valuation of $A$, the set $\ocors_v = \{ x \in A \ |\ v(x) \geq 0 \}$ is a subring, and will be called \nuovo{valuation ring} of $v$.
\end{defin}

A subring $\ocors \subset \effe$ of a field $\effe$ is said to be a \nuovo{valuation ring} for $\effe$ if $\forall x \in \effe : x \in \ocors$ or $x^{-1} \in \ocors$. If $\ocors$ is a valuation ring for $\effe$ there exists a unique equivalence class of valuations $v$ such that $\ocors = \ocors_v$ (see \cite{ZS2}). Hence the equivalence class of a valuation of a field is determined by its valuation ring.

If $v$ is a valuation of a field $\effe$, the valuation ring $\ocors_v$ is a local ring whose maximal ideal will be denoted by $m_v = \{x\ |\ v(x) > 0\}$.

\begin{teo}[(Existence of valuations)] \label{teo:esistenza posti}
Let $A$ be a ring, $k$ and $\effe$ be fields satisfying $k \subset A \subset \effe$. Let $I \subsetneq A$ be an ideal. There exists a valuation $v$ of $\effe$ such that $A \subset \ocors_v$ and $I = m_v \cap A$.

\dimo See \cite[VI\numero 4, thm.~5]{ZS2}.
\end{teo}

We will need only valuations of finite rank, thank to the following proposition.  

\begin{prop}   \label{prop:rango finito}
If $v$ is a valuation of $\effe$ over $k$, then  $\rank v$ is less or equal than the transcendence degree of $\effe$ over $k$.
\end{prop}

\subsection{Reducing the rank}  \label{subsez:abb rango}

Let $v$ be a valuation of $\effe$ of finite rank and let $\lambda \in \Lambda^v$ an elements of height $s$, i.e. $\lambda \in \Lambda^v_s \setminus \Lambda^v_{s-1}$. If we consider the elements of $\Lambda^v_{s-1}$ to be ``too small'' if compared with $\lambda$, we may take the quotient valuation:

$$\overline{v}:\effe \stackrel{v}{\freccia} \Lambda^v \cup \{+\infty\} \freccia \Lambda^v / \Lambda^v_{s-1} \cup \{+\infty\} = \Lambda^{\overline{v}} \cup \{+\infty\}$$

This way all elements of height less than $s$ become $0$, and the image of $\lambda$ in the quotient has height $1$.

We may also consider elements in $\Lambda^v \setminus \Lambda^v_s$ to be ``too big'' if compared with $\lambda$. We define a valuation of the subring $O_{v,s} = \{g \in \effe \ |\ \exists m \in \ze : v(g) \leq m \lambda\}$: 

$$\Hat{v}:O_{v,s} \ni g \freccia \left[ \left\lbrace
\begin{array}{ll}
v(g)   & \mbox{if } v(g) \in \Lambda^v_s \\
+\infty & \mbox{else}
\end{array}
\right. \right] \in \Lambda^v_s \cup \{+\infty\} = \Lambda^{\Hat{v}} \cup \{+\infty\}$$

With this valuation the elements with height greater than $s$ become $+\infty$. However the valuation $\Hat{v}$ may not be defined on the whole $\effe$, but only on the subring $O_{v,s}$. To evaluate the elements of $\effe \setminus O_{v,s}$ we should have assigned them a $-\infty$ valuation, but this is forbidden.

With this construction we may reduce the rank of a valuation, but we loose some structure on the valued set: not a field, but a ring. However if we know $\Hat{v}$ we can recover $v$, as the valuation ring $\ocors_{\Hat{v}} \subset O_{v,s}$ is a valuation ring for $\effe$, and the equivalence class of valuations corresponding to this valuation ring is the class of $v$.

Let $s$ be a fixed height, we may apply both constructions. In this way we get a rank $1$ valuation $\Hat{\overline{v}}$ of the subring $O_{v,s}$. An equivalent way for defining $\Hat{\overline{v}}$ is to take an element $\lambda \in \Lambda^v_s \setminus \Lambda^v_{s-1}$ and then put: 

$$\Hat{\overline{v}} : \effe \ni g \freccia \frac{v(g)}{\lambda} \in \erre \cup \{\pm\infty\}$$

Where $\dfrac{v(g)}{\lambda}$ is a division between elements of $\Lambda$. In this way $O_{v,s} = \effe \setminus{\Hat{\overline{v}}}^{-1}(-\infty)$.

Let $k \subset \effe$ be fields. Let $A \supset k$ be a $k$-algebra such that $\effe=\qu(A)$ (the field of fractions of $A$) and let $\famil$ be a set of generators of $A$ as a $k$-algebra.

Let $v$ be a finite rank valuation of $\effe$ over $k$ and let $s$ be the maximum height of the elements in $\famil$. Hence $A \subset O_{v,s}$, and the valuation $\Hat{\overline{v}}$ may be restricted to a rank $1$ valuation of $A$, hence it is equivalent to a real valued valuation such that $\{0\} \neq \Hat{\overline{v}}(\famil) \subset \erre$.

\begin{prop}  \label{prop:val di anelli}
Let $k,\effe,A,\famil$ be as above and let $v:A \freccia \erre \cup \{+\infty\}$ be a valuation such that $\{0\} \neq v(\famil) \subset \erre$. There exists a valuation $w$ of $\effe$ over $k$ such that $v = \Hat{\overline{w}}|_A$.

\dimo Let $P = v^{-1}(+\infty)$, a prime ideal of $A$. The valuation $v$ may be extended to the localization ring $A_P = \left\lbrace \dfrac{p}{q} \ |\ p \in A, q \in A \setminus P \right\rbrace$. We define the extension $\omega:A_P \freccia \erre \cup \{\infty\}$ by the formula $\omega\left( \frac{p}{q}\right) = v(p)-v(q)$. It is easy to see that it is well defined and that $\omega|_A = v$. The ideal $PA_P=P^e = \left\{\frac{p}{q} \ |\ p \in P, q \in A \setminus P\right\}=\{g\in A_P \ |\ \omega(g) = +\infty \}$ is the unique maximal ideal of $A_P$, hence every element $g$ with $\omega(g)\in \erre$ has an inverse in $A_P$. It follows that $\forall f,g \in A_P\setminus PA_P : \omega(f)=\omega(g) \Leftrightarrow \omega\left(\frac{f}{g}\right)=0$.

We want to define a valuation $w$ of $\cappa$ taking values in a group $\Lambda \supset \Lambda^v=\Lambda^\omega$ such that $w|_{A_P} \setminus \omega^{-1}(+\infty) = \omega$. We denote by $O$ the subring $\{f \in A_P \ |\ \omega(f) \geq 0\}$, and by $M$ the ideal $\{f \in O \ |\ \omega(f) > 0\}$. By the theorem of Existence of valuations (thm. \ref{teo:esistenza posti}) there exists a valuation ring $\ocors$ of $\effe$ with maximal ideal $m$ such that $O \subset \ocors$ and $M = O \cap m$. We denote by $w$ a valuation of $\effe$ with valuation ring $\ocors$. We want to show that it coincide with $\omega$ over $A_P$. $\forall f,g \in A_P : w(f) = w(g) \Leftrightarrow w\left(\frac{f}{g}\right)=0 \Leftrightarrow \frac{f}{g} \in \ocors \setminus m \Leftrightarrow \frac{f}{g} \in O \setminus M \Leftrightarrow \omega\left(\frac{f}{g}\right)= 0 \Leftrightarrow \omega(f)=\omega(g)$. Hence we may construct an immersion $\Lambda^v\ifreccia\Lambda^w$ such that $w|_{A_P}=\omega$.
\end{prop}

However the valuation $w$ constructed in the previous section isn't unique. 
   
\subsection{Valuating sequences}   \label{subsez:valuating seq}

Let $V$ be a complex irreducible affine variety defined over the countable field $k \subset \ci$, embedded with the classical topology. We need a way for describing valuation of $\effe = k(V)$ over $k$.

In this section we recall some definition and results from \cite{MS1} about $k$-valuating sequences. Then we introduce the concept of valuation supported out of $V'$. In the following subsections these object will be associated with the ideal points of the compactification $V'^{comp}$. 

\begin{defin}
A $k$-\nuovo{valuating} sequence is a sequence ${(x_n)}_{n \in \enne} \subset V$ such that every point $x_n$ is $k$-generic and such that:
$$ \forall f \in \effe : \exists \lim_{n \tende \infty} f(x_n) \in \ci \cup \{\infty\} = \cp^1$$
(These limits are meant in the classic topology over $\cp^1$). The condition that the points $x_n$ are $k$-generic guarantees that $f$ is well defined on them and that $f(x) \in \cappa \setminus k$.
\end{defin}

If $(x_n)$ is a valuating sequence we may define a place:
$$\plac_{(x_n)}:\effe \ni f \freccia \lim_{n \tende \infty} f(x_n) \in \cp^1$$
This function is finite over $\ocors = \{ f \in \effe \ |\ \plac_{(x_n)}(f) \neq \infty \}$. The complex function $\plac_{(x_n)} : \ocors \freccia \ci$ is a place of $\effe$ over $k$. The ring $\ocors$ is a valuation ring, and the valuation that has this ring of valuation (well defined up to equivalence) will be denoted by $v_{(x_n)}$. As the transcendence degree of $\effe$ over $k$ is finite, $v_{(x_n)}$ has finite rank.

\begin{prop}  \label{prop:esiste succ}
If $v$ is a valuation of $\effe$ over $k$, there exists a $k$-valuating sequence $(x_n)$ such that $v$ is equivalent to $v_{(x_n)}$.

\dimo \cite[corol.~1.2.3]{MS1}
\end{prop}

\begin{prop} \label{prop:logaritmi e valutazioni}
Let $(x_n) \subset V$ be a $k$-valuating sequence, and let $v = v_{(x_n)}$. Let $f,g \in \effe$ such that $v(g) \neq 0$ (i.e. $\plac_{(x_n)} g = 0$ or $\infty$). Then we may evaluate the ratio $\displaystyle \frac{v(f)}{v(g)} $ (as division between elements of $\Lambda^v$):
$$\frac{v(f)}{v(g)} =  \lim_{n \tende \infty} \frac{\log|f(x_n)|}{\log|g(x_n)|}$$

\dimo If $v(f) = 0$ the sequence $\log|f(x_n)|$ is bounded. The sequence $\log|g(x_n)|$ is not, hence both members are $0$. 

Suppose from now on that $v(f) \neq 0$.

For $h\in\effe$ we have that $v\left(\displaystyle\frac{1}{h}\right)=-v(h)$ and $\displaystyle\log\left|\frac{1}{h(x_i)}\right| = -\log|h(x_i)|$, hence we may suppose in the following that $v(f),v(g) < 0$. 

With this hypothesis we may conclude by proposition \cite[prop.~1.2.1]{MS1}. We recall the proof.

We know that $\displaystyle\lim_{n \tende \infty} \log(g(x_i)) = \infty$ and $\displaystyle\lim_{n \tende \infty} \log(f(x_i)) = \infty \Rightarrow \exists M : \forall i > M : \log|g(x_n)| > 0$ and $\log|f(x_n)| > 0$, hence $\displaystyle\frac{\log|f(x_n)|}{\log|g(x_n)|} \in (0,\infty)$.

It is enough to show only one of the two inequalities:
$$\lim_{n \tende \infty} \frac{\log|f(x_n)|}{\log|g(x_n)|} \leq \frac{v(f)}{v(g)}$$
If the second member is $0$ we have done, else we apply the same inequality exchanging the roles of $f$ and $g$, and we get the opposite inequality.

To prove the inequality we take $r,s \in \enne$ such that $\displaystyle\frac{v(f)}{v(g)} \leq \frac{r}{s}$. This implies that $s v(f) > r v(g)$ and so that $\displaystyle v\left(\frac{g^r}{f^s}\right) = r v(g) - s v(f) < 0$, hence we have $\displaystyle\log\left|\frac{g(x_i)^r}{f(x_i)^s}\right|\tende\infty$, hence $\exists M_1 : \forall i > M_1 : \displaystyle\log\left|\frac{g(x_i)^r}{f(x_i)^s}\right|\geq 0 \Rightarrow r\log|g(x_i)| \geq s\log(|f(x_i)|) \Rightarrow \frac{\log|f(x_n)|}{\log|g(x_n)|}\leq \frac{r}{s}$.
\end{prop}

\begin{defin}
A valuation $v:\effe \freccia \Lambda$ is said to be \nuovo{supported at infinity} if $k[V] \not\subset \ocors_v$.
\end{defin}

\begin{prop}
Let $v$ be a valuation of $\effe$ over $k$. $v$ is supported at infinity if and only if there exists a $k$-valuating sequence $(x_n)\subset V$ such that $v$ is equivalent to $v_{(x_n)}$, and $(x_n)$ is not contained in any compact subset of $V$.

\dimo \cite[prop.~1.2.4]{MS1}.
\end{prop}

It is easy to classify all valuations that are not supported at infinity. If $v$ is not supported at infinity, it is equivalent to a valuation $v_{(x_n)}$ where $(x_n)$ is a $k$-valuating sequence contained in a compact subset $K \subset V$. Hence we may extract a subsequence  $(x_{n_k})$ that converges to a point $x \in K$. This means that if $f \in k[V]$ there exists $\lim_{n\tende \infty} f(x_n) = f(x)$. This is true for every polynomial, hence $x_n \tende x$.

We denote by $S$ the set of all valuations of $\effe$ over $k$, up to equivalence, by $S_x$, $x \in V$ the set of all equivalence classes of valuations of the form $x_{(x_n)}$ with $x_n \tende x$, by $S_\infty$ the set of all valuation supported at infinity. With these notations we have:

$$ S = S_\infty \cup \left( \bigcup_{x \in V} S_x \right) $$

Now we fix a generating family for $V$, denoted by $\famil={(f_j)}_{j \in J} \subset k[V]$. We use the notations $V'$ and  $k[V'] = k[(f_j), (f_j^{-1})] \subset k(V)$, as usual.

\begin{defin}
A valuation $v:\effe \freccia \Lambda$ is said to be \nuovo{supported out of} $V'$ if $k[V'] \not\subset \ocors_v$, or equivalently if $\exists j : v(f_j) \neq 0$.
\end{defin}

\begin{prop}  \label{prop:supporti e successioni}
Let $v$ be a valuation of $\effe$ over $k$, $v = v_{(x_n)}$, $(x_n) \subset V$. Then $v$ is supported out of $V'$  $\Leftrightarrow \{x_n\}$ is not infinitely often contained in any compact subset of $V'$.

\dimo By definition $v$ is supported out of $V' \Leftrightarrow k[V'] \not\subset \ocors_v \Leftrightarrow$ $\exists j\in J$ such that $v(f_j) < 0$ $\Leftrightarrow \exists j : \lim f_j(x_n) = \infty$ o $\lim f_j(x_n) = 0$.

$\Rightarrow$: If $\exists j : \lim f_j(x_n) = \infty$ then $\{x_n\}$ is not infinitely often contained in a compact subset of $V$, in particular it is not infinitely often contained in a compact subset of $V'$. Else if $\exists j : \lim f_j(x_n) = 0$, the sequence $(x_n)$ may not be infinitely often contained in a compact subset of $V'$, otherwise it would have a subsequence converging to a point $x \in V'$, but then $f_j(x) \neq 0$.

$\Leftarrow$: If $(x_n)$ is not infinitely often contained in a compact subset of $V'$ we have two cases: it is infinitely often contained in a compact subset of $V$ (intersecting $V \setminus V'$) or it is not. In the former case there exists a subsequence $(x_{n_k})$ converging to a point $x \in V \setminus V'$, hence there exists $j \in J$ such that $0 = f_j(x) = \lim f_j(x_{n_k}) = \lim f_j(x_n)$. In the latter case $\exists p \in k[V] : p(x_n) \tende \infty \Rightarrow \exists j : f_j(x_n) \tende \infty$, because $\famil$ generates $k[V]$.
\end{prop}

So if we denote by $S'$ the set of equivalence classes of valuation with support out of $V'$, we can write:

$$ S' = S_\infty \cup \left( \bigcup_{x \in V \setminus V'} S_x\right)$$

\subsection{Quasi-valuating sequences}

Quasi-valuating sequences are our counterpart to Morgan and Shalen pre-valuating sequences. 

\begin{defin}
By $k$-\nuovo{quasi-valuating} sequence we mean a sequence ${(x_n)}_{n\in \enne} \subset V'$ satisfying:
\begin{enumerate}
	\item $\forall f \in k[V']: \exists \displaystyle \lim_{n \tende \infty} f(x_n) \in \ci \cup \{\infty\} = \cp^1$.
	\item $\forall f,g \in k[V'] : \displaystyle \lim_{n_\tende\infty} \log|g(x_n)| = \pm \infty \Rightarrow \exists\displaystyle \lim_{n\tende\infty} \frac{\log|f(x_n)|}{\log|g(x_n)|} \in [-\infty,\infty]$.
\end{enumerate}
\end{defin}

Note that in the previous section we used functions in $k[V']$, so they are well defined on every point of $V'$.

Let $(x_n)$ be a sequence not contained in any compact subset of $V'$. As $k[V']$ is countable by a diagonal argument we may extract a $k$-quasi-valuating subsequence not contained in any compact subset of $V'$.

\begin{prop} \label{prop:esiste succ valutante}
Let $(x_n)$ be a $k$-quasi-valuating sequence. There exists a $k$-valuating sequence $(x_n')$ such that:
\begin{enumerate}
	\item $(x_n')$ is contained in some compact subset of $V'$ if and only if $(x_n)$ is.
	\item $\forall f \in k[V'] : \displaystyle \lim_{n\tende\infty} f(x_n') = \lim_{n\tende\infty} f(x_n)$.
	\item $\forall f,g \in k[V'] : \displaystyle \lim_{n\tende\infty} \log|g(x_n)| = \pm \infty \Rightarrow \lim_{n\tende\infty} \frac{\log|f(x_n')|}{\log|g(x_n')|} = \lim_{n\tende\infty} \frac{\log|f(x_n)|}{\log|g(x_n)|}$
\end{enumerate}

\dimo As generic points are dense, and $k[V']$ is countable, we may find a sequence of generic points verifying the conditions. Then, as $\effe$ is countable, we may extract a $k$-valuating subsequence.
\end{prop}

\begin{defin}
Let $v$ be a valuation of $\effe$ over $k$ supported out of $V'$, and let $(x_n)$ be a $k$-quasi-valuating sequence not contained in any compact subset of $V'$. The valuation $v$ is said to be \nuovo{compatible} with $(x_n)$ if:
\begin{enumerate}
	\item $\displaystyle \forall f \in k[V']: \lim_{n \tende \infty} f(x_n) = \infty \Leftrightarrow v(f) < 0$.
	\item $\displaystyle \forall f \in k[V']: \lim_{n \tende \infty} f(x_n) = 0 \Leftrightarrow v(f) > 0$.
	\item $\forall f,g \in k[V'] : \displaystyle \lim_{n\tende\infty} \log|g(x_n)| = \pm \infty \Rightarrow \lim_{n\tende\infty} \frac{\log|f(x_n)|}{\log|g(x_n)|} = \frac{v(f)}{v(g)}$
\end{enumerate}
\end{defin}

\begin{prop}  \label{prop:compatibile}
For every $k$-quasi-valuating sequence $(x_n)$ not contained in any compact subset of $V'$ there exists a compatible valuation $v$ supported out of $V'$.

\dimo Let $(x_n')$ be a $k$-valuating satisfying the properties of the proposition \ref{prop:esiste succ valutante}. Let $v = v_{(x_n')}$. By proposition \ref{prop:supporti e successioni} $v$ is supported out of $V'$ and by proposition \ref{prop:logaritmi e valutazioni} is compatible with $(x_n)$.
\end{prop}

\subsection{Ideal points and valuations} \label{subsez:assocval}

We make use of the notations $S, S_\infty, S_x, S'$ defined in subsection \ref{subsez:valuating seq}.

If $v \in S$, we denote by $\Lambda^v$ its image group.The group $\Lambda^v$ has finite rank $r$, and we denote by $\Lambda^v_1 \dots \Lambda^v_r$ its convex non trivial subgroups.

We restrict our attention on values taken by elements of $\famil$.

We denote by $s$ the smaller integer such that $\forall f \in \famil : v(f) \in \Lambda^v_s$. The valuation $\overline{v}$, defined in subsection \ref{subsez:abb rango} as the composition of $v$ with the projection  $\Lambda^v \sfreccia \Lambda^v/\Lambda^v_{s-1}$ sends to $0$ all elements whose valuation has height less than $s$.

We denote by $S_0$ the set of all these valuations:

$S_0 = \{v \in S \ |\ v$ is supported out of $V'$ and $ \forall f \in \famil : v(f) \in \Lambda^v_1\}$

Let $v \in S_0$. We choose an immersion $\mu:\Lambda^v_1 \ifreccia \erre$ and we define: 

$$z = {(\mu( -v(f_j) )}_{j \in J} \in \erre^J$$

As $v$ is supported out of $V'$ we know that $z \neq 0$, hence, as $\mu$ is unique up to positive constants, it is well defined the elements $\pi(z) \in \sferic$, independently on the choose of $\mu$.

We may write $\pi(z) = U(v)$, defining a map $U:S_0 \freccia \sferic$. We want to study the image of this map, i.e. the set $U(S_0)\subset\sferic$.

Now we want to show that the set $U(S_0) \subset \sferic$ coincides with the set of ideal points $B(V')$. We need the following lemma.

\begin{lemma}          \label{lemma:limite}
Let $(x_n) \subset V'$ be a $k$-valuating sequence not contained in any compact subset of $V'$, and let $v$ be a compatible valuation of $\effe$ over $k$ supported out of $V'$. Then 

$$\displaystyle \lim_{n \tende \infty} \theta(x_n) = U(\overline{v})$$

\dimo We identify $\Lambda^{\overline{v}}_1$ with a subgroup of $\erre$, so that $U(\overline{v}) = \pi({(-\overline{v}(f_j))}_{j \in J}) \in \sferic$.

By definition $\exists j : \overline{v}(f_j) \in \Lambda^{\overline{v}}_1 \setminus {0}$. We choose $f \in \{f_j,f_j^{-1}\}$ such that $\overline{v}(f) < 0$. The existence of such an $f$ implies that $\overline{v}$ is supported out of $V'$. We have  $v(f) < 0 \Rightarrow \lim f(x_n) = \infty$ by compatibility. Hence $\exists M_1 : \forall n > M_1 : f(x_n) \neq 1$. When $n > M_1$ we define:
$$y_n = \frac{1}{\log|f(x_n)|} \Log(x_n) = {\left(\frac{\log|f_j(x_n)|}{\log|f(x_n)|}\right)}_{j \in J} \in \erre^J$$
$$y = \frac{1}{-\overline{v}(f)} U(\overline{v}) = {\left(\frac{\overline{v}(f_j)}{\overline{v}(f)} \right)}_{j \in J} \in \erre^J$$
As $v$ is supported out of $V'$, we have that $\exists M_2 > M_1 : \forall n > M_2 : y_n \neq 0$. 

We need to show that $\theta(x_n) \tende U(\overline{v})$. As $\|\theta(x_n)\| \tende \infty$, by proposition \ref{prop:sfera} it is enough to show that $\pi(\Log(x_n)) \tende U(\overline{v})$. As $\pi(\Log(x_n)) = \pi(y_n)$ and $U(\overline{v}) = \pi(y)$, it is enough to show that $y_n \tende y$. This follows by
$$ \forall j \in J : \lim_{n \tende \infty}\frac{\log|f_j(x_i)|}{\log|f(x_i)|} = \frac{v(f_j)}{v(f)} =\frac{\overline{v}(f_j)}{\overline{v}(f)} $$
The former inequality is assured by compatibility of $v$ with $(x_i)$, and the latter by definition of $\overline{v}$.
\end{lemma}

\begin{teo} \label{teo:punti ideali e valutazioni}
If $V$ is an irreducible variety, $B(V') = U(S_0)$.

\dimo $\subset$: Let $b \in B(V')$. There exists a sequence $(x_i)\subset V'$ not contained in any compact subset of $V'$, such that $\theta(x_i) \tende b$. We may extract a $k$-quasivaluating sequence. By proposition \ref{prop:compatibile} there exists a compatible valuation $v$. By previous lemma $b=U(\overline{v})$.

$\supset$: Let $v \in S_0$. By proposition \ref{prop:esiste succ} there exists a $k$-valuating sequence $(x_i)$ such that $v$ is equivalent to $v_{(x_i)}$, in particular $v$ is compatible with $(x_i)$. By previous lemma $U(v) = \lim \theta(x_i)$. By proposition \ref{prop:supporti e successioni} $(x_i)$ is not contained in any compact subset of $V'$, hence $U(v) \in B(V')$.
\end{teo}

\subsection{Restriction to real valued valuations}

The set $U(S_0)=B(V')$ may be characterized in an other way. We are interested only in the values of the elements of $\famil$, that lie in $\Lambda^{\overline{v}}_1$, so we may ignore the elements with higher values. As in subsection \ref{subsez:abb rango} we may restrict our attention to the valuation $\Hat{\overline{v}}|_{k[V']}:k[V'] \freccia \erre \cup \{\infty\}$, well defined up to a positive scalar factor.

We denote by $S_\erre$ the set of valuations $v:k[V']\freccia \erre \cup \{+\infty\}$, and we define $S_\erre' = \{v \in S_\erre \ |\ \exists j \in J : v(f_j) \neq 0 \}$. By proposition \ref{prop:val di anelli} every valuation in $S_\erre'$ may be obtained reducing the rank of a valuation in $S_0$, and reciprocally if we reduce the rank of a valuation in $S_0$ we get a valuation in $S_\erre'$ well defined up to scalar multiplication.  

We introduce the function:

$$ z:S_\erre' \ni v \freccia {(-v(f_j))}_{j \in J} \in \erre^J $$

By definition $0 \not\in z(S_\erre')$, so we may apply the projection on $\sferic$.

$$U_\erre:S_\erre' \ni v \freccia \pi(z(v)) =  \in \sferic$$

$$U_\erre(v)={[-v(f_j)]}_{j \in J}$$

With the property: $U_\erre(S_\erre') = U(S_0)$.

Now the last theorem may be stated as:

\begin{corol}
If $V$ is an irreducible variety, $U_\erre(S_\erre') = B(V')$.
\end{corol}

This result may be generalized to reducible varieties. So if $V$ is a (possibly reducible) variety, we denote by $V_i$ its irreducible components. Now we have that ${V'}^{comp} = \bigcup {V_i'}^{comp}$, $B(V') = \bigcup B(V_i')$. We introduce the notations $S_{\erre,i}' = \{v:k[V_i']\freccia \erre \cup \{+\infty\} \ |\ \exists j \in J : v(f_j) \neq 0\}$ and $S_\erre' = \{v:k[V']\freccia \erre \cup \{+\infty\} \ |\ \forall j \in J: v(f_j) \neq +\infty$ and $\exists j \in J : v(f_j) \neq 0\}$, as above. We denote by $U_{\erre,i}:S_{\erre,i}'\freccia \sferic$ and by $U_\erre: S_\erre' \freccia \sferic$ the maps defined above.

\begin{prop}   \label{prop:varrid}
If $V$ is a (possibly reducible) variety, $B(V') = U_\erre(S_\erre')$.

\dimo $\subset:$ Let $x\in B(V')$. There exists $i$ such that $x \in B(V_i') = U_{\erre,i}(S_{\erre,i}')$. Hence there exists $v:k[V_i'] \freccia \erre \cup \{+\infty\}$ such that $x = U_{\erre,i}(v)$. The map $k[V] \freccia k[V_i]$ may be extended to a surjective map $k[V'] \freccia k[V_i']$, so $k[V_i']$ is the quotient of $k[V']$ by the ideal $P_i$ extended (which will be denoted by $P_i'$). Hence $v$ may be lifted to a  valuation $w:k[V']\freccia \erre \{+\infty\}$, such that $U_\erre(w) = U_{\erre,i}(v) = x$.

$\supset:$ Let $x \in U_\erre(S_\erre')$. There exists a valuation $w \in S_\erre'$ such that $U_\erre(w) = x$. The prime ideal $P' = w^{-1}(+\infty) \subset k[V']$ is the extension of a prime ideal $P \subset k[V]$ that has empty intersection with the multiplicative set generated by $\famil$. There exists an index $i$ such that $P_i$ is contained in $P$, so by lemma \ref{lemma:quot val} there exists a valuation $v$ of $k[V_i'] = k[V'] / P_i'$ such that $x = U_{\erre,i}(v) \in B(V_i') = B(V')$
\end{prop}

The image $z(S_\erre') \subset \erre^J$ of the function $z$, united with the point $0 \in \erre^J$ is the set $CB(V')$, the cone over $B(V)'$.

\end{document}